\documentclass{amsart}
\usepackage{epsfig}

\usepackage{amssymb}
\def\leq{\leqslant}
\def\geq{\geqslant}

\newtheorem{theorem}{Theorem}
\newtheorem{proposition}[theorem]
{Proposition}
\newtheorem{lemma}[theorem]{Lemma}
\newtheorem{sublemma}[theorem]
{Sublemma}

\begin{document}

\title {Variation of the Liouville
measure of a hyperbolic surface}

\author {Francis Bonahon}
\author {Ya\c sar S\"ozen}

\address {F. Bonahon, Department
of Mathematics,  University of
Southern California, Los Angeles,
CA~90089-1113, U.S.A.}
\email{fbonahon@math.usc.edu}
\urladdr{http://math.usc.edu/\~{}%
fbonahon}

\address {Y. S\"ozen, Department
of Mathematics,  University of
Connecticut,  Storrs, CT 06269,
U.S.A.}
\email {sozen@math.uconn.edu}
\urladdr{http://www.math.uconn.edu/\~{}sozen}

\thanks{This work was partially
supported by grants DMS-9803445
and DMS-0103511 from the National
Science Foundation.}
\date{\today}

\begin{abstract} For a compact
riemannian manifold of negative
curvature,  the geodesic
foliation of its unit tangent
bundle is independent of the
negatively curved metric, up to
H\"older bicontinuous
homeomorphism. However, the
riemannian metric defines a
natural transverse measure to
this foliation, the Liouville
transverse measure, which does
depend on the metric. For a
surface $S$, we show that the map
which to a hyperbolic metric on
$S$ associates its Liouville
transverse measure is
differentiable, in an appropriate
sense. Its tangent map is valued
in the space of transverse
H\"older distributions for the
geodesic foliation.
\end{abstract}


\maketitle

One of the very basic examples of
measure preserving dynamical
system is the geodesic flow of a
riemannian manifold $S$. The
metric $m$ of $S$ has a natural
lift to a riemannian metric
$\overline m$ on the unit tangent
bundle $T^1S$, and  Liouville
observed that the volume form
defined by
$\overline m$ on $T^1S$ is
invariant under the geodesic
flow.  We want to analyze what
happens in this situation as we
vary the metric
$m$.

In general, as we modify the
metric
$m$, the topology of the geodesic
flow can dramatically change.
However, if
$S$ is compact and if $m$ has
negative curvature, a fortunate
phenomenon occurs: the \emph{geodesic foliation} of $T^1S$,
whose leaves are the orbits of
the geodesic flow, is independent
of the metric $m$. More
precisely, if $m$ and $m'$ are
two negatively curved riemannian
metrics on the compact manifold
$S$, and if $\mathcal F$ and
$\mathcal F'$ are the
corresponding geodesic foliations
of $T^1S$, there exists a
H\"older bicontinuous
homeomorphism of $T^1S$ which
sends
$\mathcal F$ to $\mathcal F'$ and
which is isotopic to the
identity; see for instance
\cite[\S 8.3]{Gro}. In
general, this homeomorphism does
not respect the parametrization of
the leaves of
$\mathcal F$ and $\mathcal F'$
provided by the geodesic flows
\cite{Ota90}.

As we move from the geodesic flow
to the geodesic foliation, there
is a one-to-one correspondence
between measures on
$T^1S$ which are invariant under
the geodesic flow and transverse
measures for the geodesic
foliation $\mathcal F$. Indeed,
any measure which is invariant
under the geodesic flow locally
is the product of some transverse
measure for $\mathcal F$ and of
the measure induced by the flow
on the leaves of $\mathcal F$. In
particular, the Liouville measure
of the metric
$m$ defines a transverse measure
$L_m$ for the geodesic foliation
$\mathcal F$ of
$m$. As indicated above, the
geodesic foliation $\mathcal F$
does not depend on the negatively
curved metric
$m$; however, as a transverse
measure for this fixed foliation,
the Liouville transverse measure
$L_m$ does depend on
$m$. 

In this paper, we restrict
attention to the case where $S$
is a compact surface and where
$m$ is a hyperbolic metric,
namely has constant curvature
$-1$. Let $\mathcal T(S)$ be the
\emph{Teichm\"uller space} of $S$,
consisting of all isotopy classes
of hyperbolic metrics on $S$, and
let
$\mathcal C(S)$ denote the space
of all transverse measures for
the geodesic foliation $\mathcal
F$ of
$T^1S$. The elements of $\mathcal
C(S)$ are called \emph{measure
geodesic currents} as they are a
special type of de Rham currents
on $T^1S$. The above construction
defines a map
\begin{equation*} L:\mathcal
T(S)\rightarrow \mathcal C(S)
\end{equation*} 
which to a
hyperbolic metric $m\in\mathcal
T(S)$ associates its Liouville
current $L_m$; indeed,
$L_m$ depends only on the isotopy
class of $m$. This map $L$ was the
focus of the article \cite{Bon88}
(see also \cite{Wol}) where, in
particular, it was shown to be a
proper embedding if
$\mathcal C(S)$ is endowed with
the weak* topology.  

The current paper is devoted to
proving that the map $L$ is
differentiable, provided that we
suitably extend its range. 

The case of the map which to
$t\in\mathbb R$ associates the
Dirac measure $\delta_t$  at the
point $t$ on $\mathbb R$ suggests
that, when considering tangent
vectors to curves in a space of
measures, one should expand one's
scope to distributions. Indeed,
in this case, the derivative
$\frac{d}{dt}\delta_t$ is the
distribution which to the
differentiable function
$\varphi$ associates
its derivative $\varphi'(t)$. 

For Liouville currents, one
encounters an additional
difficulty: Although  a
hyperbolic metric $m$ defines a
differentiable transverse
structure for the geodesic
foliation $\mathcal F$, this
transverse differentiable
structure depends on the metric
$m$. As a consequence, there is
no intrinsic general notion of
transverse distribution for
$\mathcal F$. However,
$\mathcal F$ admits an intrinsic
transverse H\"older structure,
and we can consider \emph{transverse H\"older
distributions} for $\mathcal F$.
Such a transverse H\"older
distribution assigns to  each
surface $V$ transverse to
$\mathcal F$ in
$T^1S$ a continuous linear form
on the space of compactly
supported H\"older continuous
functions
$\varphi:V\rightarrow \mathbb R$,
and this assignment is invariant
under the holonomy of $\mathcal
F$. In particular, a transverse
H\"older distribution is a
special type of transverse
distribution, and defines a
foliated de Rham current of
dimension 1 in $T^1S$, as in
\cite{Sul76}.

Let $\mathcal H(S)$ denote the
space of all transverse H\"older
distributions for the geodesic
foliation $\mathcal F$, also
called \emph{H\"older geodesic
currents}. When endowed with the
weak* topology, $\mathcal H(S)$
is an infinite dimensional
topological vector space, which
contains the space $\mathcal
C(S)$ of measure geodesic
currents.

Recall that the Teichm\"uller
space
$\mathcal T(S)$ is a
differentiable manifold of
dimension
$3\left|\chi(S)\right|$, where
$\chi(S)$ is the Euler
characteristic of $S$.

\begin{theorem} 
\label{thm:MainThm}
The Liouville map
$L:\mathcal T(S)\rightarrow
\mathcal C(S)\subset \mathcal
H(S)$ is differentiable at each
$m\in \mathcal T(S)$, in the
following sense: There exists a
linear map
$T_mL: T_m\mathcal T(S)
\rightarrow \mathcal H(S)$ such
that $T_mL(v)= \frac{d}{dt}
L_{m_t}\raise
-2pt\hbox{}_{\vert t=0}$ for every
tangent vector
$v\in T_m\mathcal T(S)$ and for
every differentiable curve
$t
\mapsto m_t\in\mathcal T(S)$,
$t\in\left]-\varepsilon,
+\varepsilon\right[$, with
$m_0=m$ and
$\frac{d}{dt}m_t\raise
-2pt\hbox{}_{\vert t=0}=v$.

In addition, the tangent map
$T_mL$ depends continuously on
$m$. 
\end{theorem}

By definition of the weak*
topology on
$\mathcal H(S)$, the statement 
$T_mL(v)= \frac{d}{dt}
L_{m_t}\raise
-2pt\hbox{}_{\vert t=0}$ means
that, for every H\"older
continuous function with compact
support
$\varphi :V \rightarrow \mathbb R$
 defined on a
surface $V$ transverse to
$\mathcal F$ in
$T^1S$, the number
$T_mL(v)(\varphi)$ associated to
$\varphi$ by the H\"older
transverse distribution $T_mL(v)$
is equal to
\begin{equation*}
T_mL(v)(\varphi)=
\frac{d}{dt}\int_V
\varphi\thinspace
dL_{m_t}\raise
-6pt\hbox{}_{\vert t=0}.
\end{equation*} 
In particular this derivative
exists.

H\"older
geodesic currents already
appeared in the paper \cite{Bon97a},
where they occurred as tangent
vectors to the space of measured
geodesic laminations on the
surface $S$; see in particular
\cite[\S\S 9--10]{Bon97a} for a few
applications of this analytic
point of view. It seems
remarkable that a notion which
was originally introduced
 because of consistency
considerations, namely because it
was the only type of transverse
distributions to $\mathcal F$
which seemed to make sense, could
actually end up being the correct
background to compute the tangent
maps of various geometric
functions. 

We prove
Theorem~\ref{thm:MainThm} by an
explicit computation using the
shearing coordinates for
Teichm\"uller space developed in
\cite{Bon96} (dual to the length
coordinates appearing in
\cite{Thu86d}). In particular, we 
give in
Theorem~\ref{thm:CompDeriv} an
explicit formula for the tangent
map $T_mL$ in terms of these
shearing coordinates. The main
part of the argument is to prove
the convergence of the series
involved in this formula.
Presumably, one could
alternatively use more classical
coordinates for $\mathcal T(S)$,
such as the Fenchel-Nielsen
coordinates. However, this seems
harder to do and the shearing
coordinates have the definite
technical advantage that one only
has to deal with shearing as
opposed to mixing lengths and
shearing.

Dragomir \v Sari\'c has recently
developed a new approach
\cite{Sar} to the results of this
article. It is somewhat more
direct, and  has the additional
advantage of extending to
(infinite dimensional)
Teichm\"uller spaces of surfaces
with infinite area, and in
particular to the universal
Teichm\"uller space. To some
extent, the occurrence of
H\"older geodesic currents in this
infinite dimensional setting is
somewhat surprising, as one could
have expected distributions with
a lower level of regularity. 

We are grateful to the referee
for a critical reading of
the manuscript and for
many insightful
suggestions. 

\section{Geodesic currents}
 Throughout the paper, $S$ will be
a compact oriented surface. The
results automatically extend to
non-orientable surfaces by
considering their orientation
coverings, and also extend to
surfaces with cusps through minor
rephrasings of the arguments.
As indicated in the introduction,
a hyperbolic metric
$m$ on $S$ defines a transverse measure
$L_m$ for the geodesic foliation
$\mathcal F$ of $T^1S$. Recall
that a transverse measure is the
assignment of a
Radon\footnote{Recall that a
\emph{Radon measure} is a measure
defined on Borel sets which
assigns finite mass to each
compact subset.} measure on each
surface
$V$ transverse to
$\mathcal F$, in a way which is
invariant under the holonomy of
$\mathcal F$. Equivalently, a
transverse measure locally is a
measure on the space of leaves of
$\mathcal F$ in local charts for
$S$, so that these measures agree
when charts overlap.

Globally, the space of leaves of
$S$ is not Hausdorff, so it is
convenient to lift the situation
to the universal covering
$\widetilde S$ of $S$. Lift the
hyperbolic metric $m$ to a
hyperbolic metric on $\widetilde
S$, which we will still denote by
$m$. The
geodesic foliation $\widetilde
{\mathcal F}$ of $T^1\widetilde
S$ defined by this metric
projects to the geodesic
foliation $\mathcal F$ of $T^1S$.
Because $\widetilde S$ is simply
connected and $m$ has negative
curvature, every geodesic is
globally length minimizing, and
it follows that $\widetilde
{\mathcal F}$ has no recurrent
leaves. As a consequence, the
space of leaves of
$\widetilde{\mathcal F}$ is
Hausdorff. This space of leaves
is also the space 
$G\bigl(\widetilde S\bigr)$ of all
oriented bi-infinite geodesics of
$\widetilde S$.

To work with $G\bigl(\widetilde
S\bigr)$, it is convenient to
consider the \emph{circle at
infinity}
$\partial_\infty\widetilde S$ of
$\widetilde S$. Recall that
$\partial_\infty\widetilde S$ is
the quotient of the set of
geodesic rays of $\widetilde S$
under the equivalence relation
which identifies two rays when
they stay at bounded distance
from each other. If we pick a
base point $O$ in $\widetilde S$,
each equivalence class has a
unique representative among the
geodesic rays issued from $O$,
which identifies
$\partial_\infty\widetilde S$ to
the unit tangent space
$T^1_O\widetilde S$. In
particular, this identification
$\partial_\infty\widetilde S\cong 
T^1_O\widetilde S$ defines a
metric on 
$\partial_\infty\widetilde S$, the
\emph{angle metric} based at $O$. 

Every geodesic of $\widetilde
S$ has two distinct end points on
the circle at infinity
$\partial_\infty\widetilde S$,
and any couple of distinct points
of
$\partial_\infty\widetilde S$ is
obtained in this way. This
identifies the space
$G\bigl(\widetilde S\bigr)$ of
oriented geodesics of
$\widetilde S$ to the open annulus
$\partial_\infty\widetilde S
\times
\partial_\infty\widetilde S
-\Delta$, where $\Delta$ denotes
the diagonal. 

A \emph{measure geodesic
current} on $S$, namely a
transverse measure for the
geodesic foliation
$\mathcal F$ in
$T^1S$, lifts to a transverse
measure for the geodesic foliation
$\widetilde{\mathcal F}$ in
$T^1\widetilde S$. It 
globally defines a measure on
the space
$G\bigl(\widetilde
S\bigr)=\partial_\infty\widetilde
S \times
\partial_\infty\widetilde S
-\Delta$ of leaves of
$\widetilde{\mathcal F}$, since
no leaf of $\widetilde{\mathcal
F}$ is recurrent. Note that this 
measure is invariant under the
natural action of the fundamental
group $\pi_1(S)$. This
identifies the space 
$\mathcal C(S)$ of measure
geodesic currents to the space of
$\pi_1(S)$--invariant measures on
$G\bigl(\widetilde S\bigr)$.   

If one
replaces the metric $m$ by
another hyperbolic metric $m'$ on
$S$ and the base point $O$ by
another point $O'\in\widetilde S$,
every
$m'$--geodesic ray issued from
$O'$ in
$\widetilde S$ is quasi-geodesic
for $m$, and therefore stays at
bounded distance from a unique
$m$--geodesic ray issued from
$O$.  This shows that the circle
at infinity
$\partial_\infty\widetilde S$ is
actually independent of the
hyperbolic metric
$m$. In addition, from the
quasi-geodesic property and from
easy curvature estimates, the two
angle metrics $d$ and $d'$ on
$\partial_\infty\widetilde S$ 
which are respectively associated
to
$m$ and the base point $O$ and to
$m'$ and the base point $O'$ are
\emph{H\"older equivalent} in the
following sense: There exists
constants
$C>0$ and
$\nu>0$ such that 
\begin{equation*}
\frac{1}{C} d(x,y)^\frac{1}{\nu}
\leq d'(x,y)
\leq C\thinspace d(x,y)^\nu
\end{equation*} for every $x$,
$y\in\partial_\infty\widetilde
S$. See also \cite[\S 7.2]{Gro}
or \cite[\S 19.1]{KatHas95}.
 In particular,
$G\bigl(\widetilde
S\bigr)=\partial_\infty\widetilde
S \times
\partial_\infty\widetilde S
-\Delta$
inherits from the angle metric of
$\partial_\infty\widetilde S$ a
natural metric well-defined up to
H\"older equivalence, and for
which the action of the
fundamental group $\pi_1(S)$ is
H\"older continuous. 

If $K$ is a compact subset of 
$G\bigl(\widetilde S\bigr)$  and
if $\nu>0$, let
$H_K^\nu\bigl(G\bigl(\widetilde
S\bigr)\bigr)$ be the set of all
H\"older continuous functions
$\varphi:G\bigl(\widetilde
S\bigr)\rightarrow \mathbb R$
with support contained in $K$ and
with H\"older exponent $\nu$.
Recall that $\varphi$ is
\emph{H\"older continuous with
H\"older exponent}
$\nu\in\left]0,1\right]$ when the
norm
\begin{equation*}
\left\|\varphi\right\|_\nu = 
\sup_{g\in G(\widetilde
S)} \left|\varphi(g)\right|+
\sup_{g,h\in G(\widetilde
S)}
\frac{\left|\varphi(g)
-\varphi(h)\right|}
{d(g,h)^\nu}
\end{equation*}
is finite. Let
$H\bigl(G\bigl(\widetilde
S\bigr)\bigr)$ be the set of
all the H\"older continuous
functions with compact support on 
$G\bigl(\widetilde S\bigr)$,
namely $H\bigl(G\bigl(\widetilde
S\bigr)\bigr)$ is the union of all
the 
$H_K^\nu\bigl(G\bigl(\widetilde
S\bigr)\bigr)$ as $K$ ranges
over all compact subsets of
$G\bigl(\widetilde S\bigr)$ and
$\nu$ ranges over all numbers in
$\left]0,1\right]$. A
\emph{H\"older geodesic current}
on $S$ is a linear map
$H\bigl(G\bigl(\widetilde
S\bigr)\bigr)\rightarrow \mathbb
R$ which is invariant under the
action of $\pi_1(S)$, and which is
continuous in the sense that its
restriction to each
$H_K^\nu\bigl(G\bigl(\widetilde
S\bigr)\bigr)$ is continuous for
the norm $\left\|~\right\|_\nu$.
Note that this notion is
independent of the choice of the
negatively curved metric $m$ on
$S$, because the distance $d$
on $G\bigl(\widetilde S\bigr)$ is
well-defined up to H\"older
equivalence.

Let $\mathcal H(S)$ denote the
space of H\"older geodesic
currents, endowed with the weak*
topology defined by the family of
semi-norms
$\alpha\mapsto
\left|\alpha(\varphi)\right|$ as
$\varphi$ ranges over all the
elements of $H\bigl(G\bigl(\widetilde
S\bigr)\bigr)$.

\section{The Liouville geodesic
current}
\label{sect:Liouv}

Let $m\in\mathcal T(S)$ be a
hyperbolic metric on $S$. Then
there is an
orientation-preserving isometry
$\Theta_m:\widetilde S\rightarrow
\mathbb H^2$ between the
universal covering $\widetilde
S$, endowed with the lift of $m$,
and the hyperbolic plane
$\mathbb H^2$. This isometry
induces H\"older bicontinuous
homeomorphisms
$\partial_\infty\widetilde S
\rightarrow
\partial_\infty\mathbb H^2$ and 
$G\bigl(\widetilde S\bigr)
\rightarrow
G\bigl(\mathbb H^2\bigr)$ which,
to avoid the proliferation of
symbols, we will also denote by
$\Theta_m$. 

Then the Liouville geodesic
current is the measure
$L_m=\Theta_m^*\bigl( L_{\mathbb
H^2} \bigr)$ on
$G\bigl(\widetilde S\bigr)$ which
is the pull back by $\Theta_m$ of
the Liouville geodesic current 
$L_{\mathbb H^2}$ on
$G\bigl(\mathbb H^2\bigr)$.

To compute $L_{\mathbb H^2}$,
consider the Poincar\'e disk
model for the hyperbolic plane
$\mathbb H^2$, namely the open
unit disk $\{z\in \mathbb C;
\left\vert z\right\vert<1\}$ with
the metric which at
$z$ is
$2/\left(1-|z|^2\right)$ times the
euclidean metric. Its circle at
infinity
$\partial_\infty
\mathbb H^2$ then is identified
to the unit circle $S^1$. If $x$,
$y\in S^1$, let $h(x,y)$ be the
oriented geodesic of $\mathbb
H^2$ going from $x$ to $y$. Then
a simple computation (see
\cite[\S A.2-3]{Bon88}) shows that
the Liouville measure $L_{\mathbb
H^2}$ on $G\bigl(\mathbb
H^2\bigr)=S^1\times S^1-\Delta$
coincides with the measure which
is locally $2dx\thinspace
dy/|x-y|^2$ at
$h(x,y)$. 

In the upper half-space model
$\{z\in\mathbb C; \mathrm{Im}(z)
>0\}$ for $\mathbb H^2$, the
circle at infinity
$\partial_\infty 
\mathbb H^2$ corresponds to
$\mathbb R\cup\{\infty\}$. Then,
the Liouville measure on
$G\bigl(\mathbb H^2\bigr)$
coincides with the measure which
is
$2dx\,dy/ \left\vert
x-y\right\vert^2$ at the geodesic
$h(x,y)$ going from $x$ to
$y\in\mathbb R$. 

The reader should beware of a few
minor differences with the
setting of \cite{Bon88}. The
Liouville geodesic current $L_m$
considered here is twice that of
\cite{Bon88}, to make it more
directly related to the volume
form of $T^1S$. Also, 
$G\bigl(\widetilde S\bigr)$
denotes here the space of
\emph{oriented} geodesics of
$\widetilde S$.

\section{The shearing coordinates
for Teichm\"uller space}
\label{sect:shearing}

We summarize here the basic
properties of the shearing
coordinates for Teichm\"uller
space which we will need. We
refer to \cite{Bon96} for details. 

These shearing coordinates are
associated to the choice of a
maximal geodesic lamination
$\lambda$ on
$S$. Recall that, for an arbitrary
metric of negative curvature on
$S$, a \emph{geodesic lamination}
is a closed subset of $S$ which
can be decomposed as a union of
disjoint complete geodesics which
have no self-intersection points.
Such a notion is actually a
topological object, independent
of the metric, in the sense that
there is a natural identification
between $m$--geodesic laminations
and $m'$--geodesic laminations
for any two negatively curved
metrics $m$ and $m'$. A geodesic
lamination is \emph{maximal} if it
is maximal for inclusion among
all geodesic laminations, which
is equivalent to the property
that the complement $S-\lambda$
consists of finitely many infinite
triangles. 

A fundamental example of a
maximal geodesic lamination is
obtained as follows. Start with a
family
$\lambda_1$ of disjoint simple
closed geodesics decomposing $S$
into pairs of pants. Each pair of
pants can be divided into two
infinite triangles by three
infinite geodesics spiralling
around some boundary components.
The union of $\lambda_1$ and of
these spiralling geodesics forms
a maximal geodesic lamination
$\lambda$.

The shearing coordinates belong
to the space $\mathcal
H(\lambda)$ of \emph{transverse
cocycles} for the maximal
geodesic lamination $\lambda$,
namely of finitely additive
signed transverse measures (valued
in
$\mathbb R$) for $\lambda$. The
vector space $\mathcal
H(\lambda)$ has finite dimension
$3|\chi(S)|$; see
\cite[Theorem~15]{Bon97b}. 
The analogy
in notation between the space
$\mathcal H(\lambda)$ of
transverse cocycles for the
geodesic lamination $\lambda$ and
the space $\mathcal H(S)$ of
H\"older geodesic currents on the
surface $S$ is not a coincidence,
because it is proved in
\cite{Bon97b} that
$\mathcal H(\lambda)$ is in a
natural way a subspace of
$\mathcal H(S)$; however, this is
irrelevant here.

The shearing
coordinates for Teichm\"uller
space define an open embedding
$\mathcal T(S)\rightarrow\mathcal
H(\lambda)$ which associates to a
hyperbolic metric $m\in\mathcal
T(S)$ its \emph{shearing
cocycle}
$\sigma_m\in\mathcal H(\lambda)$
\cite[\S2]{Bon96}. It is proved
in \cite[Theorem~A]{Bon96} that
this map defines a real analytic
diffeomorphism from $\mathcal
T(S)$ to an open subset of
$\mathcal H(\lambda)$. 

We will use a relatively
explicit description of the
inverse of the map
$m\mapsto\sigma_m$, as
constructed in \cite[\S 5]{Bon96}.
Let $m_0\in\mathcal T(S)$, let
$m$ be another hyperbolic metric
which is near $m_0$, and let
$\sigma_{m_0}$,
$\sigma_m\in\mathcal H(\lambda)$
be their respective shearing
cocycles. Then the isometries
$\Theta_{m_0}:\bigl(\widetilde
S, m_0\bigr)\rightarrow
\mathbb H^2$ and $\Theta_m:
\bigl(\widetilde
S, m\bigr)\rightarrow
\mathbb H^2$ considered in
Section~\ref{sect:Liouv} can be
chosen so that the induced map
$\Theta_m\circ\Theta_{m_0}^{-1}:
\partial_\infty\mathbb H^2
\rightarrow
\partial_\infty\mathbb H^2$ is of
the type which we now describe.
Unfortunately, this will require
some care in the
definitions as well as
relatively cumbersome notation.

Use the isometry
$\Theta_{m_0}$ to identify the
preimage $\widetilde \lambda$ of
$\lambda$ in
$\widetilde S$ to a geodesic
lamination of $\mathbb H^2$,
which we will also denote by
$\widetilde \lambda$. Pick a base
point
$O\in\widetilde S\cong\mathbb
H^2$ which is disjoint from
$\widetilde\lambda$. This defines
a partial order on the set of
components of
$\mathbb H^2-\widetilde\lambda$
which are different from the
component
$T_O$ containing $O$ where, for
two such components $T$ and $T'$, 
$T<T'$ exactly when $T$ separates
$O$ from
$T'$.

\begin{figure}[h]
\hskip 1cm
\includegraphics{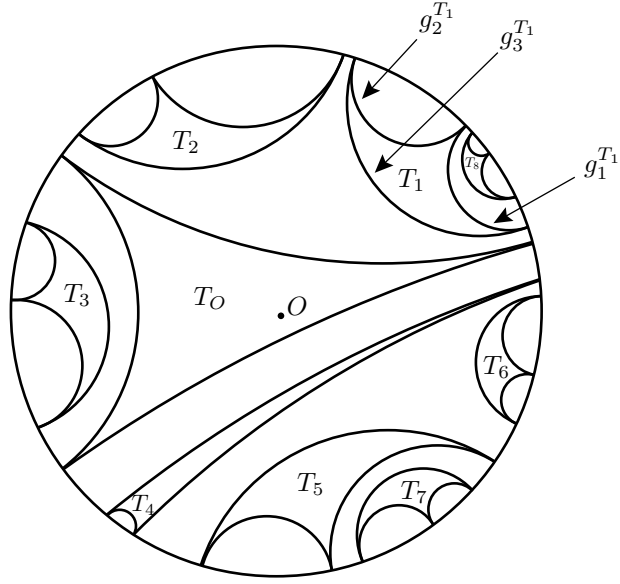}
\caption{The partial order on the
components of $\mathbb
H^2-\widetilde
\lambda$: $T_4<T_5<T_7$, $T_1$
and $T_2$ are not comparable, and
$\left\{T_1, T_3, T_6,
T_7\right\}$ is a spanning
family.}
\end{figure}

Orient the leaves of
$\widetilde\lambda$ to the left,
as seen from the base point $O$.
For such a leaf $g$ and for
$a\in\mathbb R$, let the
\emph{elementary earthquake},
or \emph{elementary shearing
map} 
$E_g^a:\mathbb H^2-g\rightarrow
\mathbb H^2-g$ coincide with the
identity on the component of
$\mathbb H^2-g$ that contains $O$
and with the hyperbolic
translation of $a$ along the
oriented geodesic $g$ on the
other component of $\mathbb
H^2-g$. Note that, although
$E_g^a$ does not continuously
extend to
$\mathbb H^2$, it has a
continuous extension to a
\emph{homeomorphism}
$\partial_\infty \mathbb H^2
\rightarrow 
\partial_\infty \mathbb H^2$,
which we also denote by $E_g^a$.

For a component $T$ of $\mathbb
H^2-\widetilde\lambda$ which is
not the component $T_O$
containing $O$, the boundary of
$T$ consists of three leaves
$g_1^T$, $g_2^T$ and $g_3^T$ of
$\widetilde\lambda$ since
$\lambda$ is a maximal geodesic
lamination. Choose the indexing
so that $g_1^T$,
$g_2^T$ and $g_3^T$ occur in this
order as one goes
counterclockwise around $\partial
T$, and so that $g_3^T$ is in the
boundary of the component of
$\mathbb H^2-T$ that contains $O$.
For $a\in\mathbb R$, define
$E_T^a=E_{g_3^T}^a\circ
E_{g_1^T}^{-a} \circ
E_{g_2^T}^{-a}$.

Set $\alpha=\sigma_m
-\sigma_{m_0}\in\mathcal
H(\lambda)$. Note that $\alpha$
also defines a transverse cocycle
for $\widetilde\lambda$. If $T$
is a component of $\mathbb
H^2-\widetilde\lambda$, let
$\alpha(T)$ be the number
associated by this transverse
cocycle to a geodesic arc
joining the base point $O$ to an
arbitrary point of $T$.

Let a \emph{spanning family} of
components of  $\mathbb
H^2-\widetilde\lambda$ be a
finite family $\mathcal U=\{U_1,
U_2, \dots, U_p\}$ of components
of $\mathbb
H^2-\widetilde\lambda$ different
from the component $T_O$
containing $O$ such that no two
$U_i$ are comparable for $<$,
namely such that no $U_i$
separates $O$ from another $U_j$.
We will write $T<\mathcal U$ when
$T<U_i$ for some $U_i\in\mathcal
U$. 

For such a spanning family
$\mathcal U$, consider the
infinite non-commutative product
\begin{equation*}
E_{\mathcal U}^\alpha =~
\overrightarrow
{\prod_{T<\mathcal U}}
E_T^{\alpha(T)}
~
\prod_{U\in\mathcal U}
E_{g_3^U}^{\alpha(U)}
\end{equation*}
defined as the limit of
\begin{equation}
\label{eqn:FiniteEarth}
E_{T_1 T_2\dots T_n}^{\alpha}=
E_{T_1}^{\alpha(T_1)}
E_{T_2}^{\alpha(T_2)}
\dots
E_{T_n}^{\alpha(T_n)}
~
\prod_{U\in\mathcal U}
E_{g_3^U}^{\alpha(U)}
\end{equation}
as the finite family $\{T_1,
T_2, \dots, T_n\}$ tends to the
set of all $T\not=T_O$ with
$T<\mathcal U$ and where, for
every $i<j$, either
$T_i<T_j$ or $T_i$ and $T_j$ are
not comparable (namely
$T_i<T_j\Rightarrow i<j$). Note
that
$E_{T_i}^{\alpha(T_i)}$ and
$E_{T_j}^{\alpha(T_j)}$ commute
when $T_i$ and $T_j$ are not
comparable, so that the order
does not matter in this case. For
the same reason, the order is
irrelevant in the product of the
(finitely many)
$E_{g_3^U}^{\alpha(U)}$ with
$U\in\mathcal U$. The arrow on
top of the first product symbol
serves as a reminder that the
terms must appear in an order
compatible with the partial order
$<$. 
The fact that the limit
exists is proved in
\cite[\S 5]{Bon96}. The set
up in that paper is slightly
different, but the restriction of
$E_{\mathcal U}^\alpha: \mathbb
H^2 - \widetilde\lambda
\rightarrow \mathbb H^2$ to a
component $T$ of $\mathbb H^2 -
\widetilde\lambda$ coincides
with what is called
$\varphi_{T_0T}$ in  \cite[\S
5]{Bon96} if $T<\mathcal U$, with
$\varphi_{T_0U_i}$ if there is an
$U_i\in \mathcal U$ such that
$U_i<T$, and with the identity
otherwise. 

We now let the spanning family
$\mathcal U$ \emph{uniformly tend
to infinity}, in the sense that
the set $\{ T; T<\mathcal U\}$
converges to the set of all
components $T\not=T_O$ of $\mathbb
H^2-\widetilde\lambda$. Then
$E_{\mathcal U}^\alpha$ converges
to a map
\begin{equation*}
E^\alpha: \mathbb
H^2-\widetilde\lambda
\rightarrow
\mathbb
H^2
\end{equation*}
which coincides with a hyperbolic
isometry on each component of $\mathbb
H^2-\widetilde\lambda$. The
convergence is immediate from the
fact that $E_{\mathcal U}^\alpha$
and $E_{\mathcal U'}^\alpha$
coincide on the component $T$
when $T<\mathcal U$ and
$T<\mathcal U'$. The map
$E^\alpha$ is the \emph{shear
map} associated to
$\widetilde\lambda$ and to the
transverse cocycle $\alpha$. 

In general, the shear map
$E^\alpha$ has no continuous
extension to $\mathbb H^2$.
However, it conjugates the action
of $\pi_1(S)$ on $\widetilde S
\cong \mathbb H^2$ to a group
$\Gamma$ of isometries of
$\mathbb H^2$ whose action is
free and properly discontinuous; 
see \cite[\S5]{Bon96}. As a
consequence, $E^\alpha$ 
continuously
extends to a homeomorphism
$E^\alpha:
\partial_\infty\mathbb H^2
\rightarrow
\partial_\infty\mathbb H^2$.
The fact
that
$\sigma_m=\sigma_{m_0} +
\alpha \in \mathcal H(\lambda)$
is the shearing cocycle of the
metric $m\in \mathcal T(S)$ means
that the isomorphism
$\pi_1(S)\cong
\Gamma$ defined by $E^\alpha$ can
be realized by an isometry between
$\left (S,m\right)$ and the
quotient
$\mathbb H^2/\Gamma$. In other
words,  the isometry $\Theta_m$
can be chosen so that the
homeomorphisms
$E^\alpha$ and
$\Theta_m\circ\Theta_{m_0}^{-1}:
\partial_\infty\mathbb H^2
\rightarrow
\partial_\infty\mathbb H^2$
coincide.

In addition, as $\mathcal U$
uniformly converges to infinity,
the homeomorphism
$E^\alpha_{\mathcal U}:
\partial_\infty\mathbb H^2
\rightarrow
\partial_\infty\mathbb H^2$
uniformly converges to
$E^\alpha:
\partial_\infty\mathbb H^2
\rightarrow
\partial_\infty\mathbb H^2$.

The homeomorphisms
$E^a_g$, $E^a_T$,
$E^\alpha_{\mathcal U}$,
$E^\alpha$ of
$\partial_\infty\mathbb H^2$
induce homeomorphisms of
$G\bigl(\mathbb
H^2\bigr)=\partial_\infty\mathbb
H^2 \times \partial_\infty\mathbb
H^2-\Delta$, which we will denote
by the same symbols.

\section{Formal computations}
\label{sect:FormalComp}

Let $t\mapsto m_t$, $t\in\left]
-\varepsilon, \varepsilon\right[$,
be a differentiable curve in the
Teichm\"uller space
$\mathcal T(S)$. We want to show
that the derivative
$\frac{d}{dt}L_{m_t}\raise
-2pt\hbox{}_{\vert t=0}$
exists as a H\"older geodesic
current. Namely, for a H\"older
continuous function
$\varphi:G\bigl(\widetilde S\bigr)
\rightarrow \mathbb R$ with
compact support, we want to show
that the derivative
\begin{equation*}
\frac{d}{dt} \iint_{G(\widetilde
S)}
\varphi\thinspace dL_{m_t}\raise
-6pt\hbox{}_{\vert t=0}
\end{equation*}
exists and depends continuously on
$\varphi$.

Let $\sigma_t\in\mathcal
H(\lambda)$ be the shearing
cocycle associated to the
hyperbolic metric, and set
$\alpha_t=\sigma_t-\sigma_0
\in\mathcal H(\lambda)$. Then,   
\begin{equation*}
 \iint_{G(\widetilde
S)}
\varphi\thinspace dL_{m_t}
=
\iint_{G(\mathbb H^2)}
\varphi\circ
\Theta_{m_t}^{-1}\thinspace
dL_{\mathbb H^2}
=
\iint_{G(\mathbb H^2)}
\varphi\circ
\Theta_{m_0}^{-1}
\circ
\left(E^{\alpha_t}\right)^{-1}
\thinspace dL_{\mathbb H^2}
\end{equation*}
where, as
in Section~\ref{sect:shearing},
the homeomorphism $E^{\alpha_t}:
G\left(\mathbb H^2\right)
\rightarrow
G\left(\mathbb H^2\right)$ is
defined as the limit of the
homeomorphisms
\begin{equation}
\label{eqn:DefPartialShear}
E_{\mathcal U}^{\alpha_t} =~
\overrightarrow
{\prod_{T<\mathcal U}}
E_T^{\alpha_t(T)}
~
\prod_{U\in\mathcal U}
E_{g_3^U}^{\alpha_t(U)}
\end{equation}
 as
the spanning family $\mathcal U$
uniformly tends to infinity. 

For a fixed
component $T$ of $\widetilde
S-\widetilde\lambda$, the number
$\sigma_t(T)$ is a linear
function of
the transverse
cocycle $\sigma_t\in\mathcal
H(\lambda)$. Since
the shearing cocycle
$\sigma_t$ depends differentiably
on the hyperbolic metric $m_t$,
it follows that the derivative
\begin{equation*}
\dot\sigma_0(T)=
\frac{d}{dt}\thinspace
\sigma_t(T)\raise
-6pt\hbox{}_{\vert t=0}=
\frac{d}{dt}\thinspace
\alpha_t(T)\raise
-6pt\hbox{}_{\vert t=0}
\end{equation*}
exists and depends only on $T$
and on the tangent vector 
$\frac{d}{dt}m_t\thinspace\raise
-2pt\hbox{}_{\vert t=0}
\in T_{m_0}\mathcal T(S)$.

If we formally differentiate the
infinite product of
Equation~\ref{eqn:DefPartialShear},
we can expect that
\begin{equation}
\label{eqn:ConjPartialShear}
\frac{d}{dt} \iint_{G(\mathbb
H^2)}
\varphi\circ
\Theta_{m_0}^{-1}
\circ
\left( E^{\alpha_t}_{\mathcal U}
\right)^{-1}
\thinspace dL_{\mathbb
H^2}\raise -6pt\hbox{}_{\vert t=0}
\kern 5pt
\raise -1pt \hbox{?}
\kern -9pt=
\sum_{T<\mathcal U}
\dot\sigma_0(T)
C_0\bigl(\varphi, T\bigr)+
\sum_{U\in\mathcal U}
\dot\sigma_0(U)
C_0\bigl(\varphi, g_3^U \bigr)
\end{equation}
where the question mark on the
$=$ sign is here to remind us
that this is only a conjectural
formula, where the contribution of
a geodesic 
$g\in G\bigl(\widetilde S\bigr)$
is
\begin{equation}
\label{eqn:DefC(g)}
C_0(\varphi, g)=
\frac{d}{da}
\iint_{G(\mathbb H^2)}
\varphi\circ
\Theta_{m_0}^{-1}
\circ
\left(E^{a}_{g}\right)^{-1}
\thinspace dL_{\mathbb
H^2}\raise -6pt\hbox{}_{\vert a=0}
\end{equation}
and where the contribution of a
component $T$ of $\widetilde S-
\widetilde\lambda$ is
\begin{equation}
\label{eqn:DefC(T)}
\begin{split}
C_0(\varphi,
T)&=\frac{d}{da}
\iint_{G(\mathbb H^2)}
\varphi\circ
\Theta_{m_0}^{-1}
\circ
\left(E^{a}_{T}\right)^{-1}
\thinspace dL_{\mathbb
H^2}\raise -6pt\hbox{}_{\vert a=0}
\\
&= C_0\bigl(\varphi, g_3^T\bigr)-
C_0\bigl(\varphi, g_1^T\bigr)-
C_0\bigl(\varphi, g_2^T\bigr)
\end{split}
\end{equation}
since $E_T^a=E_{g_3^T}^a\circ
E_{g_1^T}^{-a} \circ
E_{g_2^T}^{-a}$ by definition.

We will see that, as we let the
spanning family
$\mathcal U$ uniformly tend to
infinity, the last term 
$\sum_{U\in\mathcal U}
\dot\sigma_0(U)
C_0\bigl(\varphi, g_3^U\bigr)$ of
Equation~\ref{eqn:ConjPartialShear}
vanishes in the limit.
Therefore, the expected formula
for the derivative of $t\mapsto
L_{m_t}$ is
\begin{equation}
\label{eqn:ConjFullShear}
\begin{split}
\frac{d}{dt} \iint_{G(\widetilde
S)}
\varphi\thinspace dL_{m_t}\raise
-6pt\hbox{}_{\vert t=0}=&
\frac{d}{dt} \iint_{G(\mathbb H^2)}
\varphi\circ
\Theta_{m_0}^{-1}
\circ
\left(E^{\alpha_t}\right)^{-1}
\thinspace dL_{\mathbb
H^2}\raise -6pt\hbox{}_{\vert
t=0}\\
\kern 5pt
\raise -1pt \hbox{?}
\kern -9pt=&
\sum_{T}
\dot\sigma_0(T) C_0(\varphi,
T)
\end{split}
\end{equation}
where the sum is over all
components $T$ of $\widetilde S-
\widetilde\lambda$ that are
different from the component
$T_O$ containing the base point
$O$. 

The main technical point of our
arguments is to prove that the
infinite sums of
Equations~\ref{eqn:ConjFullShear}
and \ref{eqn:ConjPartialShear} do
converge and, unlike the above
formal computations, will
strongly use the hypothesis that
the function $\varphi$ is
H\"older continuous.

We first compute the derivative
$C_0(\varphi, g)$
corresponding to the elementary
earthquake along a single 
geodesic $g$. This will not
require much regularity for
$\varphi$.

\begin{lemma} 
\label{thm:ElemShear}
For an
oriented geodesic
$g\in G\bigl(\widetilde S\bigr)$,
let $C_0(\varphi, g)$ be defined
by Equation~\ref{eqn:DefC(g)}.
Then, for every continuous
function
$\varphi:G\bigl( \widetilde
S\bigr) \rightarrow \mathbb R$
with compact support,
\begin{equation}
\label{eqn:ElemShear}
C_0(\varphi, g) = \iint
_{G(\widetilde S)}
\varphi(h)\thinspace\cos\theta(g,h)
\thinspace dL_{m_0} (h)
\end{equation} 
where
$\theta(g,h)\in\left[0,\pi\right[$
denotes the angle from $g$ to
$h$ at their intersection point,
measured counterclockwise for the
metric
$m_0$ and without taking the
orientation of the geodesics
into account, and where by
convention
$\cos\theta(g,h)=0$ 
if $g$ and $h$ do not meet. 
\end{lemma}
\begin{proof} Before going any 
further, we will introduce a
restriction on the function
$\varphi$ which will somewhat
simplify the exposition.
The Liouville measures
$L_m$ are invariant under the
involution $r$ of
$G\bigl(\widetilde S\bigr)$
defined by reversing the
orientation of the geodesics. As a
consequence, the integral of
$\varphi$ for $L_m$ is equal to
the integral of the symmetrized
function
$(\varphi+\varphi\circ r)/2$.
The same applies to the integral
of Equation~\ref{eqn:ElemShear}.
We can therefore restrict
attention to those functions
$\varphi$ which are
\emph{balanced} in the sense that
they are invariant under the
involution $r$.

Choose the isometric
identification
$\Theta_{m_0}:\bigl(\widetilde S,
m_0\bigr)\rightarrow \mathbb H^2$
so that, in the upper half-space
model for $\mathbb H^2$, the
oriented geodesic $g$ goes from
the point $0$ to the point
$\infty$. Then, on
$\partial_\infty \mathbb H^2 =
\mathbb R\cup\{\infty\}$, the
homeomorphism $E^a_g$ is such
that  $E^a_g(x)=\mathrm{e}^ax$ for
 $x\geq0$ and $E^a_g(y)=y$ for
$y\leq0$. Since only those
geodesics which cross $g$
contribute to the derivative, it
follows that
\begin{equation*}
\begin{split}
C_0(\varphi, g) &= \frac{d}{da}
\iint_{G(\mathbb H^2)}
\varphi\circ
\Theta_{m_0}^{-1}
\circ
\left(E^{a}_{g}\right)^{-1}
\thinspace dL_{\mathbb
H^2}\raise -6pt\hbox{}_{\vert
a=0}\\
&= 4\frac{d}{da}
\int_{-\infty}^0
\int_0^{+\infty}
\varphi\circ
\Theta_{m_0}^{-1}
\left(h\left(
\mathrm{e}^{-a} x, y
\right)\right)
\thinspace 
\frac{dx\,dy}{\left(
x-y\right)^2}
\raise
-6pt\hbox{\phantom{m}}_{\vert
a=0}\\
&= 4\frac{d}{da}
\int_{-\infty}^0
\int_0^{+\infty}
\varphi\circ
\Theta_{m_0}^{-1}
\left(h\left(
 x, y
\right)\right)
\thinspace 
\frac{\mathrm{e}^{a}
dx\,dy}{\left(
\mathrm{e}^{a}x-y\right)^2}
\raise
-6pt\hbox{\phantom{m}}_{\vert
a=0}\\
&= 4
\int_{-\infty}^0
\int_0^{+\infty}
\varphi\circ
\Theta_{m_0}^{-1}
\left(h\left(
 x, y
\right)\right)
\frac{-x-y}{x-y}
\frac{dx\,dy}{\left(
x-y\right)^2}
\\
\end{split}
\end{equation*}
where $h(x,y)$ denotes the
geodesic joining $x$ to $y$. Here
the factor $4$ has two origins:
one factor $2$ comes from the
expression of the Liouville
geodesic current in
Section~\ref{sect:Liouv}; another
$2$ comes from the fact that we
have to consider the two possible
orientations for each geodesic
crossing $g$, and that each
orientation contributes the same
amount by our assumption that
$\varphi$ is balanced. An
elementary computation in
hyperbolic geometry shows that
the term
$\frac{-x-y}{x-y}$ is equal to
the cosine of the angle from the
geodesic $g$ joining $0$ to
$\infty$ to the geodesic $h(x,y)$
joining $x$ to $y$. This
concludes the proof.
\end{proof}

\section{The main estimate}
\label{sect:MainEst}

Throughout this section, $S$ is
endowed with the hyperbolic
metric $m_0$, and $\lambda$ is
realized by an $m_0$--geodesic
lamination. We pick a base point
$O\in\widetilde S$ which is
disjoint from the preimage
$\widetilde\lambda$ of $\lambda$.
Endow the circle at infinity
$\partial_\infty\widetilde S$
with the angle metric as seen
from the point
$O$, and the space 
$G\bigl(\widetilde S\bigr)
=\partial_\infty\widetilde S
\times \partial_\infty\widetilde S
-\Delta$ of oriented
geodesics of
$\widetilde S$ with the
corresponding product metric.
For a H\"older continuous function
$\varphi:G\bigl(\widetilde
S\bigr)\rightarrow\mathbb R$, its
H\"older norm
$\left\Vert\varphi\right\Vert_\nu$
will be computed with respect to
this metric on 
$G\bigl(\widetilde S\bigr)$.

Consider a component $T$ of
$\widetilde S-\widetilde\lambda$
which is not the component $T_O$
containing the base point $O$.
Let $D_T$ denote the distance
from $O$ to $T$, namely the
distance from $O$ to the point
$x_T$ of the geodesic
$g_3^T$ that is closest to
$O$. Let $y_T$ be the orthogonal
projection to
$g_3^T$ of the vertex of the
triangle
$T$ that separates $g_1^T$ from
$g_2^T$. Finally, let $u_T\in
\mathbb R$ be the signed distance
from
$x_T$ to
$y_T$ in the oriented geodesic
$g_3^T$, namely $\left\vert u_T
\right\vert$ is equal to
the distance from
$x_T$ to $y_T$ and
$u_T<0$ exactly when $y_T$ is to
the right of $x_T$ as seen from
$O$.

\begin{figure}[h]
\includegraphics{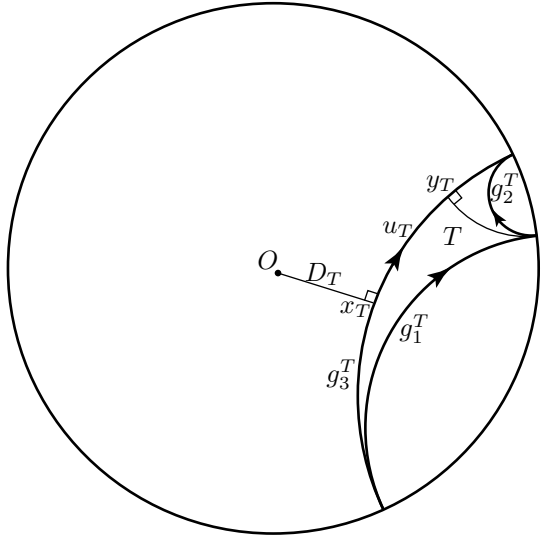}
\caption{The quantities $D_T$ and
$u_T$.}
\end{figure}

The key ingredient of the proof
of Theorem~\ref{thm:MainThm} is
the following estimate for the
quantity $C_0(\varphi, T)$
defined by
Equation~\ref{eqn:DefC(T)}.

We write $u\prec v$ when the ratio
$\left\lvert u/v\right\rvert$ is
bounded above by a constant, and
we write
$u\asymp v$ when both
$u\prec v$ and $v\prec u$.

\begin{proposition}
\label{thm:EstC(T)}
Let
$\varphi:G(\widetilde
S)\rightarrow\mathbb R$ be a
H\"older continuous function with
H\"older exponent $\nu$, and such
that every geodesic $g\in
G(\widetilde S)$ in the support
of $\varphi$ meets the ball
$B(O,R)\subset\widetilde S$.
Then,
\begin{equation*}
 C_0(\varphi, T)\prec
\left\Vert\varphi\right\Vert_\nu
\mathrm{e}^{-(1+\nu)D_T}
\mathrm{e}^{-\left\vert
u_T\right\vert}
\end{equation*}
 where the constant hidden in the
symbol $\prec$ depends only on
the radius $R$ and on a number
$r>0$ such that
$\widetilde\lambda$ is disjoint
from
$B(O,r)$.
\end{proposition}

We can rephrase the estimate of
Proposition~\ref{thm:EstC(T)} as
follows. The hyperbolic triangle
$T$ has a natural \emph{center}
$O_T$, defined for instance as
the point that is at equal
distance from the three sides of
$T$.  

\begin{lemma}
\label{thm:EstO_T}
\begin{equation*}
\bigl\vert
d(O,O_T)-D_T - \left\vert
u_T\right\vert
\bigr\vert
\prec 1
\end{equation*}
\end{lemma}
\begin{proof}[Proof of
Lemma~\ref{thm:EstO_T}]
By elementary hyperbolic
geometry, the distance 
$d(O,O_T)$ is within
$\frac12\log 3$ of $d(O,y_T)$
which, by consideration of the
right angled hyperbolic triangle
$Ox_Ty_T$, is
itself within $\log
2$ of $D_T+\left\vert
u_T\right\vert$.
\end{proof}

The estimate of
Proposition~\ref{thm:EstC(T)} can
therefore be rewritten as
\begin{equation*}
 C_0(\varphi,
T)\prec
\left\Vert\varphi\right\Vert_\nu
\mathrm{e}^{-\nu d(O,T)}
\mathrm{e}^{-d(O,O_T)}.
\end{equation*}

We now begin the proof of
Proposition~\ref{thm:EstC(T)},
which will occupy most of this
section. We split the argument in
two cases, according to whether
the distance
$D_T$ from $O$ to $T$ is larger
than
$2R$ or not.
\begin{proof}[Proof of
Proposition~\ref{thm:EstC(T)}
when
$D_T\geq 2R$]  Choose the
isometry
$\Theta_{m_0}$ between
$\widetilde S$, endowed with the
metric $m_0$, and the Poincar\'e
disk model for the hyperbolic
space $\mathbb H^2$, bounded by
the unit circle
$S^1$ in $\mathbb C$, so that
$\Theta_{m_0}$ sends the base
point
$O\in\widetilde S$  to
the center $0\in\mathbb
H^2\subset \mathbb C$ of the
disk. We use $\Theta_{m_0}$ to
identify $\widetilde S$ to
$\mathbb H^2$ so that, in this
identification $\widetilde
S\cong\mathbb H^2$, the map
$\Theta_{m_0}$ is now the
identity. We will consequently use
the same symbol to represent an
object in
$\widetilde S$ and the
corresponding object in $\mathbb
H^2$.

As in the proof of
Lemma~\ref{thm:ElemShear}, we can
assume without loss of
generality that $\varphi$ is
balanced, namely invariant under
the involution $r$ of
$G\bigl(\widetilde S\bigr)$
defined by reversing the
orientation of geodesics. 

Every geodesic $h$ which
contributes to the integrals in
$C_0(\varphi, T)$ must be in the
support of
$\varphi$ and must cross at least
one of the geodesics
$g_i^T$. Because of our
assumption that
$D_T\geq 2R$ and by choice of $R$,
such an
$h$ will necessarily cross
$g_3^T$. The vertices of the
ideal triangle
$T$ decompose the circle
$\partial_\infty \widetilde S\cong
S^1$ into three intervals $I_1$,
$I_2$ and
$I_3$, where $I_i$ has the same
end points as the geodesic
$g_i^T$. If we combine the
definition of $C_0(\varphi, T)$ in
Equation~\ref{eqn:DefC(T)}, the
computation of
Lemma~\ref{thm:ElemShear} for
$C_0(\varphi, g)$, the expression
for the Liouville geodesic current
$L_{\mathbb H^2}$ given in
Section~\ref{sect:Liouv}, and 
the fact that
$\varphi$ is balanced, it follows
that
\begin{multline*}
 C_0(\varphi,
T)=4\int_{I_3}\int_{I_1\cup I_2}
\varphi(h(x,y))\thinspace
\left(\cos\theta(g_3^T,h(x,y))
\right.\\
-
\left.\cos\theta(g_1^T,h(x,y))-
\cos\theta(g_2^T,h(x,y))\right)
\frac{dx\thinspace dy}{\left\lvert
x-y\right\rvert^2}
\end{multline*}
 where $h(x,y)$ denotes the
geodesic with end points $x$,
$y\in S^1$.  It therefore
suffices to bound the quantity
\begin{multline*} I(y)=
2\int_{I_1\cup I_2}
\varphi(h(x,y))\thinspace
\left(\cos\theta(g_3^T,h(x,y))
\right.\\
-\left.
\cos\theta(g_1^T,h(x,y))-
\cos\theta(g_2^T,h(x,y))\right)
\frac{dx}{\left\lvert
x-y\right\rvert^2}
\end{multline*}
 for every $y\in I_3$. (The
factor $2$ is here to avoid
dragging cumbersome constants in
the computations below.)

We now switch to the upper-half
space model for $\mathbb H^2$,
bounded by the real line $\mathbb
R$ in $\mathbb C$, in such a way
that the base point $O$ now
corresponds to the point
$i\in\mathbb H^2\subset\mathbb C$
and in such a way that the point
$y$ is sent to $\infty$. The
measure
$dx/\left\lvert x-y\right\rvert^2$
on
$S^1$ is then sent to the measure
$\frac12 dx$ on $\mathbb
R$. In this context, this gives
the expression
\begin{multline*}
 I(\infty)= \int_{a}^{b}
\varphi(h(x,\infty))\thinspace
\left(
\cos\theta(g_3^T,h(x,\infty))
\right.\\
-\left.
\cos\theta(g_1^T,h(x,\infty))-
\cos\theta(g_2^T,h(x,\infty))\right)
\thinspace dx
\end{multline*}
 where, as before, $h(x,\infty)$
is the geodesic going from $x$ to
$y=\infty$, namely the vertical
euclidean line passing through
$x$, and where $a$ and
$b\in\mathbb R$ correspond to the
end points of
$I_3$, with $a<b$.

Let $c\in\mathbb R$ correspond to
the third vertex of $T$. If we use
elementary euclidean geometry to
compute
$\cos\theta(g_k^T,h(x,\infty))$
and if we write 
$\psi(x)=\varphi(h(x,\infty))$,
then 
\begin{multline*}
 I(\infty) =
\int_{a}^{b}
\psi(x)  \frac{2x-a-b}{b-a} dx
\\
-
\int_{a}^{c}
\psi(x)  \frac{2x-a-c}{c-a}
 dx-
\int_{c}^{b}
\psi(x) \frac {2x-c-b}{b-c} dx.
\end{multline*}
 Note that
\begin{equation*}
\int_{a}^{b} \frac{2x-a-b}{b-a}
dx=
\int_{a}^{c} \frac{2x-a-c}{c-a}
 dx=
\int_{c}^{b} \frac{2x-c-b}{b-c}
dx=0.
\end{equation*}  Therefore,
\begin{equation*}
\begin{split}
I(\infty) &=
\int_{a}^{b}
\eta(x)  \frac{2x-a-b}{b-a} dx-
\int_{a}^{c}
\eta(x)  \frac{2x-a-c}{c-a}
 dx-
\int_{c}^{b}
\eta(x)  \frac{2x-c-b}{b-c} dx\\
&=-
\int_{a}^{c}
\eta(x)
\frac{2(x-a)(b-c)}{(b-a)(c-a)} dx
+ \int_{c}^{b}
\eta(x)
\frac{2(b-x)(c-a)}{(b-a)(b-c)} dx
\end{split}  
\end{equation*}
if we set
$\eta(x)=\psi(x)-\psi(c)$.
Integrating in $x$, we conclude
that
\begin{equation}
\label{eqn:EstI(infty)}
\left\vert I(\infty)\right\vert 
\leq 
2 \frac{(c-a)(b-c)}{b-a}
\max_{x \in
[a,b]}\left\vert\eta(x)\right\vert
\end{equation} 
We now use the fact that 
$\varphi:G(\widetilde
S)\rightarrow
\mathbb R$ is H\"older continuous
with respect to the angle metric
on $S^1$, with H\"older exponent
$\nu\in\left]0,1\right]$ and with
H\"older norm
$\left\Vert\varphi\right\Vert_\nu$.
As we switch from the disk model
to the upper half-space model for
$\mathbb H^2$ by an isometry
sending $0$ to $i$, the angle
metric on
$S^1$ becomes the metric $
\frac{2}{1+x^2}\thinspace dx$ on
$\mathbb R$. It follows that the
function
$\eta(x)=\psi(x)-\psi(c)=
\varphi(h(x,\infty))-
\varphi((h(c,\infty))$ is H\"older
continuous with respect to the
usual metric $dx$ of $\mathbb R$,
with H\"older exponent $\nu$
and with H\"older norm
$\left\Vert\eta\right\Vert_\nu\leq
2\left\Vert\varphi\right\Vert_\nu$.

In particular, because
$\eta(c)=0$, we have that
$\left\vert\eta(x)\right\vert
\leq
2\left\Vert\varphi\right\Vert_\nu
\left\vert x-c\right\vert^\nu$ for
every $x$. It follows that 
$\max_{x \in
[a,b]}\left\vert\eta(x)\right\vert$
is bounded by 
$2\left\Vert\varphi\right\Vert_\nu
\left\vert b-a\right\vert^\nu$. 
Combining  with
Equation~\ref{eqn:EstI(infty)},
we obtain that 
\begin{equation}
\label{eqn:Est(infty)}
\left\vert I(\infty)\right\vert
\leq
4\left\Vert\varphi\right\Vert_\nu
\frac{(c-a)(b-c)}{b-a}
\left\vert b-a\right\vert^\nu.
\end{equation} 

We now estimate $D_T$ and $u_T$ in
terms of $a$, $b$ and $c$. 

\begin{lemma}
\label{thm:Estu_tAndD_t} 
If
$I(\infty)\not=0$ and if $T$ is at
distance at least
$2R$ from the base point $O$,
then:
\begin{enumerate}
\item[(i)] 
$\mathrm{e}^{-D_T}\asymp b-a;$
\item[(ii)]  $c-a
\asymp b-a$ and
$\mathrm{e}^{-\left\vert
u_T\right\vert}\asymp
\frac{b-c}{b-a}$ when
$u_T\geq0;$
\item[(iii)]  $b-c
\asymp b-a$ and
$\mathrm{e}^{-\left\vert
u_T\right\vert}\asymp
\frac{c-a}{b-a}$ when $u_T\leq0$.
\end{enumerate}
\end{lemma}

\begin{proof}[Proof of
Lemma~\ref{thm:Estu_tAndD_t}] On
the geodesic $g_3^T$, with end
points $a$ and $b$, we already
have two preferred points, namely
the orthogonal projection $x_T$
of $O=i$ and the orthogonal
projection $y_T$ of the vertex
$c$. Recall that
$D_T$ is the distance from $O$ to
$x_T$, and that
$u_T$ is the signed distance from
$x_T$ to $y_T$ in the oriented
geodesic $g_3^T$. Consider a
third point, namely the
orthogonal projection
$z_T$ of the point
$\infty$ to $g_3^T$. In
euclidean terms,
$z_T$ is the point of the euclidean semi-circle
$g_3^T\subset
\mathbb H^2
\subset \mathbb C$ that has the
largest imaginary part, namely
$z_T=\frac{a+b}{2}+\frac{b-a}{2}i$.

\begin{figure}[h]
\includegraphics{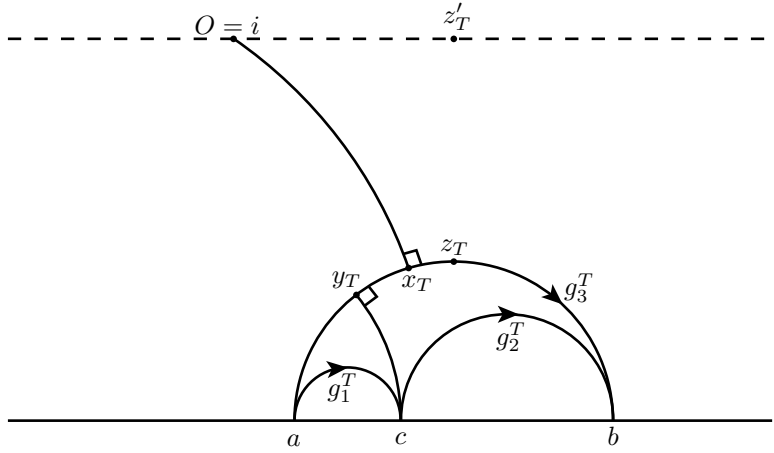}
\caption{The proof of
Lemma~\ref{thm:Estu_tAndD_t},
with $u_T<0$.}
\end{figure}

The points $x\in\mathbb R$ such
that the geodesic joining $x$ to
$\infty$ meets the disk $B(O,R)$
form an interval
$[-R', R']$, where $R'$ depends
only on $R$. If
$I(\infty)\not=0$, then the
intervals
$[a,b]$ and $[-R', R']$ have a
non-trivial intersection. Since,
in addition,
$g^T_3$ is disjoint from the
interior of the ball $B(O,2R)$,
we conclude that
$a$ and
$b$ stay in a closed interval
$[-R'', R'']$ where, again,
$R''$ depends only on $R$. 

Consider the horocycle centered at
$\infty$ passing through $O$,
namely the euclidean horizontal
line passing through $O=i$. Let
$z_T'$ be the point of this
horocycle which lies on the same
vertical line as $z_T$, namely
$z_T'=\frac{a+b}{2}+i$. The piece
of horocycle joining $O$ to
$z_T'$ has
  length $\left\vert
\frac{a+b}{2}\right\vert
\leq R''$. It follows that the
distance from $O$ to
$z_T'$ is bounded by $R''$. Also,
the orthogonal projection of
$z_T'$ to
$g_3^T$ is equal to $z_T$, the
orthogonal projection of $O$ to
$g_3^T$ is equal to $x_T$, and the
orthogonal projection map is
distance non-increasing; it
follows that the distance from
$x_T$ to $z_T$ is also bounded by
$R''$. 

By definition, $D_T$ is equal to
the distance from $O$ to $x_T$.
From the above observations and
the triangle inequality, it
follows that $D_T$ is within
$2R''$ of the distance from
$z_T'=\frac{a+b}{2}+i$ to
$z_T=\frac{a+b}{2}+\frac{b-a}{2}i$.
Since this distance is equal to
$\left\vert
\log\frac{b-a}{2}\right\vert$,
 and since $\frac{b-a}{2}\leq
R''$ this proves that
$\mathrm{e}^{-D_T}\asymp b-a$.

Similarly, $u_T$ is equal to the
signed distance from the point
$x_T$ to $y_T$ in $g_3^T$, and is
therefore within $R''$ of the
signed distance
$u_T'$ from $z_T$ to $y_T$. 
Consider the  isometry
$\Theta$ of $\mathbb H^2$ sending
$a$, $b$, $\infty$ to
$\infty$, $0$, $1$, respectively.
Then
$\Theta$ sends $c$ to
$-\mathrm{e}^{-u_T'}$. From the
invariance of cross-products under
$\Theta$, we conclude that
\begin{equation*}
\frac{(-\mathrm{e}^{-u_T'}-0)(1-\infty)}
{(-\mathrm{e}^{-u_T'}-\infty)(1-0)}=
\frac
{(c-b)(\infty-a)}{(c-a)(\infty-b)}
\end{equation*} namely
$\mathrm{e}^{-u_T'}=\frac{b-c}{c-a}$.
It follows that
\begin{equation*}
\mathrm{e}^{-u_T}\asymp
\mathrm{e}^{-u_T'}=\frac{b-c}{c-a}.
\end{equation*}

If $u_T\geq0$, then $u_T'$ admits
a uniform lower bound and
$\mathrm{e}^{-u_T'}$ admits a
uniform upper bound, so that 
$b-c\prec c-a$. It follows that
$c-a
\asymp b-a$. This implies
\begin{equation*}
\mathrm{e}^{-\left\vert
u_T\right\vert}=
\mathrm{e}^{-u_T}\asymp
\frac{b-c}{c-a}
\asymp \frac{b-c}{b-a}.
\end{equation*}

Similarly, when $u_T\leq0$, then
$b-c\asymp b-a$ and 
\begin{equation*}
\mathrm{e}^{-\left\vert
u_T\right\vert}=
\mathrm{e}^{u_T}\asymp
\mathrm{e}^{u_T'}=\frac{c-a}{b-c}
\asymp \frac{c-a}{b-a}.
\end{equation*}

This concludes the proof of
Lemma~\ref{thm:Estu_tAndD_t}.

\end{proof}

We now return to our estimate 
\begin{equation*}
\left\vert I(\infty)\right\vert
\leq 4\left\Vert\varphi\right
\Vert_\nu
\frac{(c-a)(b-c)}{b-a}
\left\vert b-a\right\vert^\nu
\end{equation*}
from
Equation~\ref{eqn:Est(infty)}.

When $u_T\geq0$,
Lemma~\ref{thm:Estu_tAndD_t} and
the fact that $c-a\leq b-a$ imply
\begin{equation*}
\begin{split}
\left\vert I(\infty)\right\vert
&\leq
4\left\Vert\varphi\right\Vert_\nu
(b-c)
\left\vert b-a\right\vert^\nu\\
&\prec 
\left\Vert\varphi\right\Vert_\nu
\mathrm{e}^{-D_T}
\mathrm{e}^{-\left\vert
u_T\right\vert}
\mathrm{e}^{-\nu D_T}.
\end{split}
\end{equation*}
 The same estimate
holds when $u_T\leq0$ by a
symmetric argument. 

If we go back to the beginning of
the proof, this shows that 
\begin{equation*}
I(y) \prec
\left\Vert\varphi\right\Vert_\nu
\mathrm{e}^{-(1+\nu)D_T}
\mathrm{e}^{-\left\vert
u_T\right\vert}
\end{equation*}
 for every $y\in
I_3$. Then,
\begin{equation*}
 C_0(\varphi,
T) = 2 
\int_{I_3} I(y)
\thinspace dy
\prec
\left\Vert\varphi\right\Vert_\nu
\mathrm{e}^{-(1+\nu)D_T}
\mathrm{e}^{-\left\vert
u_T\right\vert}.
\end{equation*}

 This concludes the proof of
Proposition~\ref{thm:EstC(T)}
when $T$ is at distance at least
$2R$ from the base point $O$.
\end{proof}

\begin{proof}[Proof of
Proposition~\ref{thm:EstC(T)}
when $D_T\leq 2R$] Note that, in
this case, we only need to show
that 
$ C_0(\varphi,
T)
\prec
\left\Vert\varphi\right\Vert_\nu
\mathrm{e}^{-\left\vert
u_T \right\vert}$. Also, we again
can assume that
$\varphi$ is balanced, without
loss of generality.

In this case, we have an
additional term because there can
be geodesics
$h$ which contribute to
$C_0(\varphi, T)$ but which do
not meet the geodesic
$g_3^T$. These additional
geodesics must cross both $g_1^T$
and $g_2^T$, and therefore must
have one end point in the
interval $I_1$ and the other one
in $I_2$. So, $C_0(\varphi, T)$
splits as a sum
$C_0(\varphi, T)=C_1(\varphi,
T)+C_2(\varphi, T)$ with
\begin{multline*} 
C_1(\varphi,
T)=4\int_{I_3}\int_{I_1\cup I_2}
\varphi(h(x,y))\thinspace
\left(\cos\theta(g_3^T,h(x,y))
\right.\\
- \left.
\cos\theta(g_1^T,h(x,y))-
\cos\theta(g_2^T,h(x,y))\right)
\frac{dx\thinspace dy}{\left\lvert
x-y\right\rvert^2}
\end{multline*}
 and
\begin{multline*}
 C_2(\varphi,
T)=4\int_{I_1}\int_{I_2}
\varphi(h(x,y))\thinspace
\left(-\cos\theta(g_1^T,h(x,y))
\right.\\
- \left.
\cos\theta(g_2^T,h(x,y))\right)
\frac{dx\thinspace dy}{\left\lvert
x-y\right\rvert^2}
\end{multline*}
using the fact that $\varphi$ is
balanced. 

As in the first case, to bound the
first integral $C_1(\varphi, T)$,
it suffices to bound  the quantity
\begin{multline*} I(y)=
2\int_{I_1\cup I_2}
\varphi(h(x,y))\thinspace
\left(\cos\theta(g_3^T,h(x,y))
\right.\\
- \left.
\cos\theta(g_1^T,h(x,y))-
\cos\theta(g_2^T,h(x,y))\right)
\frac{dx}{\left\lvert
x-y\right\rvert^2}
\end{multline*}
 for every $y\in I_3$. 

For future reference, let us note
the following easy fact.
\begin{lemma}
\label{thm:I(y)bounded}
\begin{equation*}
 I(y)
\prec 
\left\Vert \varphi
\right\Vert_\nu.
\end{equation*}
\end{lemma}

\begin{proof}[Proof of
Lemma~\ref{thm:I(y)bounded}]
If
$x$ contributes to $I(y)$, then
$h(x,y)$ meets the disk $B(O,R)$,
so that
$\left\lvert x-y\right\rvert$ is
bounded away from 0. \end{proof}

Let us switch again to the upper
half-space model for $\mathbb
H^2$, in such a way that $y$
corresponds to
$\infty$ and $O$ corresponds to
$i$. Then $I(y)$ becomes
$I(\infty)$ with
\begin{equation}
\label{eqn:FormI(infty)}
\begin{split}
 I(\infty) &=
\int_{a}^{b}
\psi(x)  \frac{2x-a-b}{b-a} dx-
\int_{a}^{c}
\psi(x)  \frac{2x-a-c}{c-a}
 dx-
\int_{c}^{b}
\psi(x)  \frac{2x-c-b}{b-c} dx\\
&=-
\int_{a}^{c}
\psi(x)
\frac{2(x-a)(b-c)}{(b-a)(c-a)} dx
+
\int_{c}^{b}
\psi(x)
\frac{2(b-x)(c-a)}{(b-a)(b-c)} dx
\end{split}
\end{equation}
 with
$\psi(x)=\varphi(h(x,\infty))$ (no
need to introduce a function
$\eta$ now). If we again note that
$\left\vert\psi(x)
\right\vert\leq 
\left\Vert\varphi\right\Vert_\nu$
for every $x$ and if we integrate
in $x$, it follows that 
\begin{equation}
\label{EstI(infty)bis}
\left\vert I(\infty)\right\vert
\leq
2\left\Vert\varphi\right\Vert_\nu
\frac{(c-a)(b-c)}{b-a}.
\end{equation}

Let us relate this estimate to
$\mathrm{e}^{-\left\vert
u_T\right\vert}$.

\begin{lemma} 
\label{thm:Estu_tAndD_tBis}
If $T$
is at distance at most $2R$ and at
least
$r>0$ from the base point
$O$, 
\begin{enumerate}
\item[(i)]
$\mathrm{e}^{-\left\vert
u_T\right\vert}\asymp
\frac{b-c}{1+bc}$ when
$u_T\geq0;$
\item[(ii)]
$\mathrm{e}^{-\left\vert
u_T\right\vert}\asymp
\frac{c-a}{1+ac}$ when
$u_T\leq0$.
\end{enumerate}
\end{lemma}
\begin{proof}[Proof of
Lemma~\ref{thm:Estu_tAndD_tBis}]
Consider the case where $u_T\leq
0$. The other case is similar.

Let $\rho$ be the hyperbolic
rotation around the point
$O=i$ that sends $a$ to 0. Namely,
$\rho(z)=(z-a)/(az+1)$. The
distance from $O$ to the geodesic
$\rho(g_3^T)$ is equal to the
distance from $O$ to
$g_3^T$ and, by hypothesis, is
therefore bounded between the two
positive constants $r$ and
$2R$. Since the end points of
$\rho(g_3^T)$ are
$\rho(a)=0$ and $\rho(b)$, it
follows that there are two
constants
$0<R_1<R_2<\infty$ such that
$\rho(b)\in[R_1, R_2]$. Namely
$\rho(b)\asymp 1$. 

Since $\rho(a)$ and $\rho(b)$ stay
bounded, we can apply to
$\rho(T)$ the arguments of the
proof of
Lemma~\ref{thm:Estu_tAndD_t}. 
Note that
$\left\vert u_T\right\vert$ is
equal to the distance from
$\rho(x_T)$ to
$\rho(y_T)$. On the euclidean
semi-circle
$\rho(g_3^T)$, let $z_T$ be the
point that has the largest
imaginary part and let $u_T'$ be
the signed distance from
$\rho(y_T)$ to $z_T$ in
$\rho(g_3^T)$. Then, as in the
proof of
Lemma~\ref{thm:Estu_tAndD_t},
$u_T$ is within a uniformly
bounded distance of $u_T'$ and
\begin{equation*}
\mathrm{e}^{-\left\vert
u_T\right\vert}=
\mathrm{e}^{u_T}\asymp
\mathrm{e}^{u_T'}
=\frac{\rho(c)-\rho(a)}{\rho(b)-\rho(c)}
=\frac{\rho(c)}{\rho(b)-\rho(c)}
\end{equation*} since
$\rho(a)=0$. Because
$u_T\leq0$, it follows that
$\rho(c)\prec\rho(b)-\rho(c)$, and
therefore that
$\rho(b)-\rho(c)\asymp
\rho(b)\asymp 1$. This proves
\begin{equation*}
\mathrm{e}^{-\left\vert
u_T\right\vert}\asymp
\frac{\rho(c)}{\rho(b)-\rho(c)}
\asymp \rho(c) =\frac{c-a}{1+ac}.
\end{equation*}

This concludes the proof of
Lemma~\ref{thm:Estu_tAndD_tBis}
when $u_T\leq0$. The other case
is similar, using the hyperbolic
rotation around $O$ that sends
$b$ to 0.

\end{proof}

We now return to $I(\infty)$. We
 want to prove that
$ I(\infty)\prec
\left\Vert\varphi\right\Vert_\nu 
\mathrm{e}^{-\left\vert
u_T\right\vert}$.

The support of
$\psi$ is contained in the
interval
$[-R', R']$ consisting of those
$x$ such that the geodesic
joining $x$ to
$\infty$ meets the disk $B(O,R)$.
Therefore, if
$I(\infty)\not=0$, the interval
$[a,b]$ must meet $[-R', R']$.
Also, if $[a,b]$ contains $[-R',
R']$, its length $b-a$ cannot be
too large, as
$O=i$ would otherwise be on the
wrong side of the oriented
geodesic $g_3^T$ going from $a$
to $b$. Therefore, there exists
$R''>R'$, depending only on
$R$, such that at least one of
$a$,
$b$ belongs to the interval
$[-R'', R'']$. 

Let us restrict attention, without
loss of generality, to the case
where
$a$ is in $[-R'', R'']$. We now
split the argument in two cases,
according to whether $u_T\leq 0$
or
$u_T\geq0$. 

First consider the case where
$u_T\leq 0$. Since $I(\infty)$ is
bounded by
Lemma~\ref{thm:I(y)bounded}, we
can restrict attention to the
case where
$\left\vert u_T\right\vert$ is
uniformly large, namely where, by
Lemma~\ref{thm:Estu_tAndD_tBis}, 
$w=(c-a)/(1+ac)>0$ is uniformly
small. Note that
$c-a=w(1+a^2)/(1-aw)$ so that
$c-a\asymp w$ if
$0<w<1/(2R'')$ (since $-R''\leq
a\leq R''$). By
Equation~\ref{EstI(infty)bis} and
Lemma~\ref{thm:Estu_tAndD_tBis},
it then follows that

\begin{equation*}
\left\vert I(\infty)\right\vert
\leq
2\left\Vert\varphi\right\Vert_\nu
\frac{(c-a)(b-c)}{b-a}
\leq
2\left\Vert\varphi\right\Vert_\nu
(c-a)
\asymp
\left\Vert\varphi\right\Vert_\nu w
\asymp
\left\Vert\varphi\right\Vert_\nu 
\mathrm{e}^{-\left\vert
u_T\right\vert}
\end{equation*}
 which concludes the proof that 
$\left\vert
I(\infty)\right\vert\prec
\left\Vert\varphi\right\Vert_\nu 
\mathrm{e}^{-\left\vert
u_T\right\vert}$  in this case.

Now consider the case where
$u_T\geq0$. Again, by
Lemmas~\ref{thm:I(y)bounded} and
\ref{thm:Estu_tAndD_tBis}, we can
restrict attention to the case 
where
$\left\vert u_T\right\vert$ is
uniformly large, namely where 
$w'=(b-c)/(1+bc)>0$ is uniformly
small.

If $b$ is uniformly bounded, say
$\left\vert b\right\vert\leq
2R''$, the same argument as above
gives that $b-c\asymp w'$ if
$0<w'<1/(4R'')$, so that
\begin{equation*}
\left\vert I(\infty)\right\vert
\leq
2\left\Vert\varphi\right\Vert_\nu
\frac{(c-a)(b-c)}{b-a}
\leq
2\left\Vert\varphi\right\Vert_\nu
(b-c)
\asymp
\left\Vert\varphi\right\Vert_\nu
w'
\asymp
\left\Vert\varphi\right\Vert_\nu 
\mathrm{e}^{-\left\vert
u_T\right\vert}.
\end{equation*}

In the other case where
$\left\vert b\right\vert\geq
2R''$, then
$b\geq 2R''$ since $b>a\geq
-R''$. Note that
\begin{equation*}
c=(b-w')/(1+bw')\geq
(2R''-w')/(1+2R''w')
\geq R''> R'
\end{equation*}
 if
$w'<\min\{\frac{1}{4R''},
\frac{R''}{2}\}$. (For the first
inequality, note that the
function $b\mapsto
(b-w')/(1+bw')$ is monotone
increasing). We now need to go
back to the original expression
for
$I(\infty)$ in
Equation~\ref{eqn:FormI(infty)}.
If we remember that the support of
$\psi$ is contained in $[-R',
R']$,
\begin{equation*}
\begin{split}
\left|I(\infty)\right|
&=\left|- \int_{a}^{c}
\psi(x)
\frac{2(x-a)(b-c)}{(b-a)(c-a)} dx
+
\int_{c}^{b}
\psi(x)
\frac{2(b-x)(c-a)}{(b-a)(b-c)}
dx\right|\\ 
&=\left| \int_{a}^{R'}
\psi(x)
\frac{2(x-a)(b-c)}{(b-a)(c-a)}
dx\right|\\ 
&\leq
\left\Vert\varphi\right\Vert_\nu
(R'-a)^2
\frac{b-c}{(c-a)(b-a)}\\ 
&\leq
\left\Vert\varphi\right\Vert_\nu
(R'-a)^2 
\frac{b-c}{1+bc}
\frac{1+bc}{(c-a)(b-a)}\\ 
&\prec 
\left\Vert\varphi\right\Vert_\nu
\mathrm{e}^{-\left\vert
u_T\right\vert}
\end{split}
\end{equation*}
 by
Lemma~\ref{thm:Estu_tAndD_tBis}
and using the fact that 
$\frac{1+bc}{(c-a)(b-a)}$ is
uniformly bounded when $a\leq
R'<R''\leq c<b$.

This proves that 
$I(\infty)\prec
\left\Vert\varphi\right\Vert_\nu 
\mathrm{e}^{-\left\vert
u_T\right\vert}$ in all cases.

Therefore  
$I(y)\prec
\left\Vert\varphi\right\Vert_\nu
\mathrm{e}^{-\left\vert
u_T\right\vert}$ for every $y\in
I_3$, and
\begin{equation*}
 C_1(\varphi,
T) =2 
\int_{I_3} I(y)
\thinspace dy
\prec
\left\Vert\varphi\right\Vert_\nu
\mathrm{e}^{-\left\vert
u_T\right\vert}.
\end{equation*}

We now consider the other term
\begin{equation*}
\begin{split}
C_2(\varphi, T)
 &=
4\int_{I_1}\int_{I_2}
\varphi(h(x,y))\thinspace
\left(-\cos\theta(g_1^T,h(x,y))-
\cos\theta(g_2^T,h(x,y))\right)
\frac{dx\thinspace dy}{\left\lvert
x-y\right\rvert^2}\\
&\prec 
\left\Vert\varphi\right\Vert_\nu
\mathrm{length}(I_1)~
\mathrm{length}(I_2)
\end{split}
\end{equation*}
 since $\left\lvert
x-y\right\rvert$ is bounded away
from 0 when the geodesic
$h(x, y)$ belongs to the support
of
$\varphi$.

\begin{lemma}
\label{thm:LengthEst}
~ 
\begin{enumerate}
\item[(i)]
$\mathrm{length}(I_1)\asymp
\mathrm{e}^{-\left\vert
u_T\right\vert}
$ and
$\mathrm{length}(I_2)\asymp 1$ 
when
$u_T\leq 0;$
\item[(ii)]
$\mathrm{length}(I_2)\asymp
\mathrm{e}^{-\left\vert
u_T\right\vert}
$ and
$\mathrm{length}(I_1)\asymp 1$ 
when
$u_T\geq 0$.
\end{enumerate}
\end{lemma}

\begin{proof}[Proof of
Lemma~\ref{thm:LengthEst}] By
symmetry, we can restrict
attention to the case where
$u_T\leq 0$ without loss of
generality.

The length of $I_1\subset S^1$ is
equal to the angle under which
$I_1$ is seen from $O$. Since the
projection
$x_T$ of $O$ to $g_3^T$ stays in
the ball $B(O,2R)$, this angle
is comparable to the angle
$\theta$ under which
$I_1$ is seen from $x_T$, namely 
$\mathrm{length}(I_1)\asymp
\theta$. Considering the
right-angled hyperbolic 
triangle whose vertices are
$x_T$, the vertex $c$ of $T$ that
separates
$I_1$ from $I_2$, and the
projection
$y_T$ of $c$ to $g_3^T$, 
elementary hyperbolic
trigonometry gives
\begin{equation*}
 \mathrm{length}(I_1)\asymp
\sin\theta =\frac{1}{\cosh
\left\vert u_T\right\vert}
\asymp
\mathrm{e}^{-\left\vert
u_T\right\vert}
\end{equation*}
 since $\left\vert u_T\right\vert$
is equal to the hyperbolic
distance from $x_T$ to
$y_T$ by definition.

Because the distance from $O$ to
$g_3^T$ is bounded,
$\mathrm{length}(I_1)+
\mathrm{length}(I_2)\asymp 1$, so
that
$\mathrm{length}(I_2)\asymp 1$.

This concludes the proof of
Lemma~\ref{thm:LengthEst}. 
\end{proof}

From Lemma~\ref{thm:LengthEst}
and our prior estimate, we obtain
\begin{equation*}
C_2(\varphi, T)
\prec 
\left\Vert\varphi\right\Vert_\nu
\mathrm{length}(I_1)\,
\mathrm{length}(I_2)
\prec
\left\Vert\varphi\right\Vert_\nu
\mathrm{e}^{-\left\vert
u_T\right\vert}.
\end{equation*}
Since we already showed that 
$ C_1(\varphi, T)
\prec
\left\Vert\varphi\right\Vert_\nu
\mathrm{e}^{-\left\vert
u_T\right\vert}$, this proves
\begin{equation*}
C_0(\varphi,
T)=
C_1(\varphi, T)+C_2(\varphi,
T)
\prec
\left\Vert\varphi\right\Vert_\nu
\mathrm{e}^{-\left\vert
u_T\right\vert},
\end{equation*}
which concludes the proof of
Proposition~\ref{thm:EstC(T)}
when $D_T\leq 2R$, and
therefore in all cases.
\end{proof}

We will need a
refinement of
Proposition~\ref{thm:EstC(T)},
which is automatically given by
the above proof. Indeed, it did
not use the full force of the
hypothesis that
$\varphi$ is H\"older continuous.
The H\"older continuity of
$\varphi$ did not come up at all
in the case where $D_T\leq 2R$,
where only the
$\mathrm L^\infty$ norm
$\sup_{g\in G(\widetilde
S)}\left\vert \varphi(g)
\right\vert \leq \left\Vert
\varphi \right\Vert_\nu$ was
used. In the case where $D_T\geq
2R$, the only occurrence of the
H\"older continuity was just
above
Equation~\ref{eqn:Est(infty)},
and more precisely in the
statement that
\begin{equation*}
\left\vert \eta(x)\right\vert
= \left\vert \psi(x)-\psi(c)
\right\vert =\left\vert
\varphi(h(x, \infty)) -
\varphi (h(c, \infty))\right\vert
\leq 2\left\Vert
\varphi\right\Vert_\nu
\left\vert b-a \right\vert^\nu
\end{equation*}
for every $x \in \left[a,
b\right]$, with the notation used
there. Because of the
estimate that $b-a \asymp \mathrm
e^{-D_T}$ from
Lemma~\ref{thm:Estu_tAndD_t}, we
can replace this inequality by
the stronger (up to a constant)
statement that
\begin{equation*}
\left\vert \eta(x)\right\vert
\prec
\left\Vert \varphi
\right\Vert_{g_3^T, \nu} \mathrm
e^{-\nu D_T}
\end{equation*}
where, for a geodesic $g\in
G\bigl( \widetilde S \bigr)$
which is disjoint from the base
point $O$, 
$\left\Vert
\varphi
\right\Vert_{g, \nu} \
\prec\left\Vert \varphi
\right\Vert_\nu$ is defined by 
\begin{equation}
\label{Eqn:DefTNuNorm}
\left\Vert \varphi
\right\Vert_{g, \nu} =
\sup_{h\in G(\widetilde
S)} \left|\varphi(h)\right|+
\sup_{h,h'}
\frac{\left|\varphi(h)
-\varphi(h')\right|}
{ \mathrm
e^{-\nu D_g}}
\end{equation}
where the second supremum is
taken over all geodesics $h$,
$h'$ which both
cross the geodesic $g$ and
have the same end point on the
side of $g$ facing the base point
$O$, and where $D_g$ is the
distance from $O$ to $g$. 

Then the proof of
Proposition~\ref{thm:EstC(T)}
automatically gives:

\begin{lemma}
\label{thm:EstC(T)bis}
Let
$\varphi:G(\widetilde
S)\rightarrow\mathbb R$ be a
H\"older continuous function with
H\"older exponent $\nu$, and such
that every geodesic $g\in
G(\widetilde S)$ in the support
of $\varphi$ meets the ball
$B(O,R)\subset\widetilde S$.
Then,
\begin{equation*}
C_0(\varphi,
T)\prec
\left\Vert\varphi\right\Vert_{g_3^T,
\nu}
\mathrm{e}^{-(1+\nu)D_T}
\mathrm{e}^{-\left\vert
u_T\right\vert}
\end{equation*}
 where the constant hidden in the
symbol $\prec$ depends only on
the radius $R$ and on a number
$r>0$ such that
$\widetilde\lambda$ is disjoint
from
$B(O,r)$. \qed
\end{lemma}

We will need a similar estimate
for the quantity
\begin{equation*}
C_0(\varphi, g)= \frac{d}{da}
\iint_{G(\mathbb H^2)}
\varphi\circ
\Theta_{m_0}^{-1}
\circ
\left(E^{a}_{g}\right)^{-1}
\thinspace dL_{\mathbb
H^2}\raise -6pt\hbox{}_{\vert a=0}
\end{equation*}
associated to the geodesic $g\in
G\bigl(\widetilde S\bigr)$ as in
Equation~\ref{eqn:DefC(g)}.
Again, let
$D_g$ be the distance from
$g$ to the base point $O$,
measured in the metric $m_0$.
\begin{proposition}
\label{thm:EstC(g)}
Let
$\varphi:G(\widetilde
S)\rightarrow\mathbb R$ be a
H\"older continuous function with
H\"older exponent $\nu$, and such
that every geodesic $g\in
G(\widetilde S)$ in the support
of $\varphi$ meets the ball
$B(O,R)\subset\widetilde S$.
Then,
\begin{equation*}
\left\vert C_0(\varphi,
g)\right\vert\prec
\left\Vert\varphi\right\Vert_{g,
\nu}
\mathrm{e}^{-(1+\nu)D_g}\prec
\left\Vert\varphi\right\Vert_\nu
\mathrm{e}^{-(1+\nu)D_g}
\end{equation*}
 where the constant hidden in the
symbol $\prec$ depends only on
the radius $R$.
\end{proposition}
\begin{proof} The proof is very
similar to, and much simpler than,
the first case of the proof of
Proposition~\ref{thm:EstC(T)}. We
consequently omit it. 
\end{proof}

\section {Convergence of infinite
sums} 
This section is devoted to
proving that the conjectural
formula of
Equation~\ref{eqn:ConjFullShear}
really makes sense. 

The arguments
of Section~\ref{sect:MainEst} were
exclusively based in the
universal cover $\widetilde
S\cong \mathbb H^2$. In this
section, we will strongly use the
action of
$\pi_1(S)$ and the fact that the
quotient space $S$ of $\widetilde
S$ under this action is compact. 

\begin{proposition} 
\label{thm:SumCv}
The sum
\begin{equation*}
\sum_{T} \dot\sigma_0(T)
C_0(\varphi, T)
\end{equation*}
 is convergent, and its absolute
value is bounded by
$c\left\Vert
\varphi\right\Vert_\nu$ for a
constant $c$ depending only on the
H\"older exponent $\nu$ and on a
compact subset containing the
support of $\varphi$ (as well as
on the metric $m_0\in \mathcal
T(S)$ and on the transverse
cocycle
$\dot\sigma_0
\in\mathcal H(\lambda)$). 
\end{proposition}

\begin{proof} By
Proposition~\ref{thm:EstC(T)}, 
it suffices to show that 
$\sum_{T} 
\left|\dot\sigma_0(T)\right|
\mathrm{e}^{-(1+\nu)D_T}
\mathrm{e}^{-\left|u_T\right|}$
converges.

 Let $\Phi\subset S$ be a
train track carrying
$\lambda$. Recall that a
\emph{(fattened) train track}
$\Phi$ on the surface $S$ is a
family of finitely many
 `long' rectangles $e_1$, \dots,
$e_n$  which are foliated
 by arcs parallel to the `short'
sides and which meet only along
arcs (possibly reduced to a
point) contained in their short
sides. In addition, a train track
$\Phi$ must satisfy the
following: 
\begin{enumerate}
 \item[(i)]  each point of the
`short' side of a rectangle also
belongs to  another rectangle,
and each component of the union
of the short sides of all
rectangles is an arc,  as opposed
to a closed curve;
 \item[(ii)]  note that the
closure $\overline{S-\Phi}$ of
the complement $S-\Phi $ has a
certain number of `spikes',
corresponding to the points where
at least 3 rectangles meet; we
require that no component of
$\overline{S-\Phi}$ is a disc
with 0, 1 or 2 spikes or an
annulus with no spike. 
\end{enumerate}
The rectangles are called the
\emph{edges} of $\Phi$. The
foliations of the edges of $\Phi$
induce a foliation of $\Phi$,
whose leaves are the \emph{ties}
of the train track. The finitely
many ties where several edges
meet are the \emph{switches} of
the train track
$\Phi$. A tie which is not a
switch is \emph{generic}. 
The geodesic lamination $\lambda
$ is \emph{carried} by 
the train track $\Phi$ if it is
contained in the interior of 
$\Phi$ and if its leaves are
transverse to the ties of $\Phi$.
There are several
constructions which provide a
train track $\Phi$ carrying
$\lambda$; see for instance
\cite{PenHar}\cite{CanEpsGre}.  
 
Let $\widetilde\Phi$ be the
preimage of $\Phi$ in the
universal covering
$\widetilde S$.
Choose an
orientation for the ties of each
edge of $\Phi$ (with no
requirement that these
orientations match at the
switches), and lift this
orientation to all the edges of
$\widetilde\Phi$. For each edge
$\widetilde e$ of
$\widetilde\Phi$, let
$T_{\widetilde e}$ be the
component of
$\widetilde S-\widetilde\lambda$
that contains the negative end
points of the ties of $\widetilde
e$, for this choice of
orientation. 

Recall that, for a component
$T\not=T_O$ of $\widetilde
S-\widetilde\lambda$, $x_T$ is the
point of the closure of $T$ that
is closest to the base point $O$.
To simplify the exposition, we
can arrange by a small
perturbation of $\Phi$ that no
$x_T$ in contained in a switch of
$\widetilde\Phi$. As a
consequence, each $x_T$ is
contained in a unique edge
$\widetilde e$ of
$\widetilde\Phi$. Then,
\begin{equation}
\label{equn:11.1}
\begin{split}
\sum_{T} \left|\dot\sigma_0(T)
C_0(\varphi, T)\right| &\prec
\left\Vert\varphi\right\Vert_\nu
\sum_{T} 
\left|\dot\sigma_0(T)\right|
\mathrm{e}^{-(1+\nu)D_T}
\mathrm{e}^{-\left|u_T\right|}\\
&\prec
\left\Vert\varphi\right\Vert_\nu
\sum_{\widetilde e}
\sum_{x_T\in \widetilde e}
\left|\dot\sigma_0(T)\right|
\mathrm{e}^{-(1+\nu)D_T}
\mathrm{e}^{-\left|u_T\right|}\\
&\prec
\left\Vert\varphi\right\Vert_\nu
\sum_{\widetilde e}
\sum_{x_T\in \widetilde e}
\left|\dot\sigma_0(T)
-\dot\sigma_0(T_{\widetilde
e})\right|
\mathrm{e}^{-(1+\nu)D_T}
\mathrm{e}^{-\left|u_T\right|}\\
&\qquad\qquad\qquad+
\left\Vert\varphi\right\Vert_\nu
\sum_{\widetilde e}
\sum_{x_T\in \widetilde e}
\left|\dot\sigma_0(T_{\widetilde
e})\right|
\mathrm{e}^{-(1+\nu)D_T}
\mathrm{e}^{-\left|u_T\right|}
\end{split}
\end{equation}
 where the sums
$\sum_{\widetilde e}$ are over all
edges of
$\widetilde\Phi$ and where the
sums
$\sum_{x_T\in \tilde e}$ are over
those components $T$ of
$\widetilde S-\widetilde\lambda$
for which $x_T$ is in $\widetilde
e$.

The distance $D_T$ is also the
distance from $x_T$ to the
base point $O$. It follows that,
if
$x_T\in \widetilde e$, then 
\begin{equation}
\label{equn:11.2}
\mathrm{e}^{-D_T}
\asymp
\mathrm{e}^{-D_{\tilde e}}
\end{equation}
 where $D_{\tilde e}$ is the
distance from $O$ to $\widetilde
e$ and where the constants hidden
in the symbol
$\asymp$ depend only on the
diameters of the edges
$\widetilde e$, which are
uniformly bounded by equivariance
under the action of
$\pi_1(S)$.

If $x_T\in \widetilde e$,
consider the \emph{divergence
radius} $r(T)$ of the component
$T\cap\widetilde e$ of
$\widetilde e-\widetilde \lambda$,
defined as follows. If
$T\cap\widetilde e$ is one of the
two components of $\widetilde
e-\widetilde
\lambda$ that touch the boundary
$\partial\widetilde\Phi$, then
$r(T)=1$. Otherwise,
$T\cap\widetilde e$ is delimited
by two geodesic arcs contained in
two of the three leaves $g_i^T$,
$i\in\{1, 2, 3\}$,  of
$\widetilde\lambda$ bounding $T$.
These two leaves are asymptotic
on one side of $\widetilde e$. On
the other side, they cross the
same succession of $r$ edges of
$\widetilde \Phi$ before taking
diverging routes at some switch of
$\widetilde \Phi$. Then, $r(T)$ is
defined as this number $r$ of
edges which the two leaves cross
together before diverging.

A consequence of the
combinatorics of transverse
cocycles, proved in
\cite[Lemma~6]{Bon97b}, is that 
\begin{equation}
\label{equn:11.3}
\left|\dot\sigma_0(T)-
\dot\sigma_0(T_{\tilde e})\right|
\prec
\left\| \dot\sigma_0\right\|_\Phi
r(T)\prec r(T)
\end{equation}
where $\left\|
\dot\sigma_0\right\|_\Phi
=\max\left\{ 
\left\vert \dot\sigma_0 (e)
\right\vert;\, e
\mathrm{~edge~of~}
\Phi \right\}$ and where we
absorb this quantity $\left\|
\dot\sigma_0\right\|_\Phi$ in the
constant of the second $\prec$. 

The complement $S-\lambda$ has
finitely many components; see for
instance
\cite{CanEpsGre}\cite{PenHar}. As
a consequence, there are only
finitely many
$T$ modulo the action of
$\pi_1(S)$. It follows that, when
we consider the point $y_T$
projection to
$g_3^T$ of the third vertex of
$T$, the distance from $y_T$ to
$T\cap\partial\widetilde\Phi$ is
uniformly bounded. From the
definition of $r(T)$ and the
definition of $u_T$ as the signed
distance from $x_T$ to
$y_T$, we conclude that
$r(T)\asymp \left\vert
u_T\right\vert +1$, where the
constants in the symbol
$\asymp$ depend only on the
(suitably defined) lengths of the
edges of $\Phi$. It follows that
there is an $a>0$ such that
\begin{equation}
\label{equn:11.4}
\mathrm{e}^{-\left|u_T\right|}
\prec
\mathrm{e}^{-ar(T)}.
\end{equation}

Finally, using a compact
fundamental domain for the
action of $\pi_1(S)$ on
$\widetilde S$, one can go from
$O$ to
$T_{\tilde e}$ by a path which
crosses at most $\prec D_{\tilde
e}$ edges of
$\widetilde\Phi$, or more
precisely which is made up of
$\prec D_{\tilde e}$ ties and of
arcs in
$\widetilde S-\widetilde\Phi$. It
follows that 
\begin{equation}
\label{equn:11.5}
\left|
\dot\sigma_0(T_{\tilde e})\right|
\prec
\left\| \dot\sigma_0\right\|_\Phi
D_{\tilde e}\prec D_{\tilde e}.
\end{equation}  

Combining  the estimates of
Equations~\ref
{equn:11.1}--\ref{equn:11.5}, we
obtain

\begin{equation*}
\begin{split}
\sum_{T} \left|\dot\sigma_0(T)
C_0(\varphi, T)\right| &\prec 
\left\Vert\varphi\right\Vert_\nu
\sum_{\widetilde e}
\sum_{x_T\in \widetilde e}
\left|\dot\sigma_0(T)
-\dot\sigma_0(T_{\widetilde
e})\right|
\mathrm{e}^{-(1+\nu)D_T}
\mathrm{e}^{-\left|u_T\right|}\\
&\qquad\qquad\qquad+
\left\Vert\varphi\right\Vert_\nu
\sum_{\widetilde e}
\sum_{x_T\in \widetilde e}
\left|\dot\sigma_0(T_{\widetilde
e})\right|
\mathrm{e}^{-(1+\nu)D_T}
\mathrm{e}^{-\left|u_T\right|}\\
&\prec
\left\Vert\varphi\right\Vert_\nu
\sum_{\widetilde
e}
\sum_{x_T\in \widetilde e} r(T)
\mathrm{e}^{-(1+\nu)D_{\tilde e}}
\mathrm{e}^{-ar(T)}\\
&\qquad\qquad\qquad+
\left\Vert\varphi\right\Vert_\nu
\sum_{\widetilde e}
\sum_{x_T\in \widetilde e}
D_{\tilde e}
\mathrm{e}^{-(1+\nu)D_{\tilde e}}
\mathrm{e}^{-ar(T)}\\ &\prec 
\left\Vert\varphi\right\Vert_\nu
\sum_{\widetilde e}
\sum_{x_T\in \widetilde e}
\mathrm{e}^{-(1+\nu')D_{\tilde e}}
\mathrm{e}^{-a'r(T)}
\end{split}
\end{equation*}
 for arbitrary $\nu'<\nu$ and
$a'<a$.

Given $\widetilde e$ and an
integer
$r\geq 1$, the number of $T$
meeting $\widetilde e$ (or,
equivalently of components
$T\cap\widetilde e$ of
$\widetilde e
-\widetilde\lambda$) with
$r(T)=r$ is uniformly bounded by
a number depending only on the
topology of
$S$; see \cite[Lemma~4]{Bon96}.
It follows that 
\begin{equation*}
\sum_{x_T\in \widetilde e}
\mathrm{e}^{-a'r(T)}
\prec
\sum_{r=1}^\infty 
\mathrm{e}^{-a'r}
\prec 1
\end{equation*}
 and therefore that 
\begin{equation*}
\sum_{T} \left|\dot\sigma_0(T)
C_0(\varphi, T)\right|
\prec 
\left\Vert\varphi\right\Vert_\nu
\sum_{\widetilde e}
\mathrm{e}^{-(1+\nu')D_{\tilde
e}}.
\end{equation*}

Now, lift the edges $e_i$, $i=1$,
\dots, $n$, of
$\Phi$ to edges $\widetilde e_i$
of
$\widetilde\Phi$, so that every
edge
$\widetilde e$ is of the form
$\gamma\widetilde e_i$ for some
$\gamma\in\pi_1(S)$ and for one of
these finitely many $\widetilde
e_i$. If 
$\widetilde e=\gamma\widetilde
e_i$, then
$\mathrm{e}^{-D_{\tilde e}}
\asymp
\mathrm{e}^{-d(O,\gamma O)}$, so
that 
\begin{equation*}
\sum_{T}\left| \dot\sigma_0(T)
C_0(\varphi, T)\right|
\prec 
\left\Vert\varphi\right\Vert_\nu
\sum_{\widetilde e}
\mathrm{e}^{-(1+\nu')D_{\tilde e}}
\prec 
\left\Vert\varphi\right\Vert_\nu
\sum_{i=1}^n
\sum_{\gamma\in\pi_1(S)}
\mathrm{e}^{-(1+\nu')d(O,\gamma
O)} .
\end{equation*}
 By convergence of the Poincar\'e
series
\begin{equation*}
\sum_{\gamma\in\pi_1(S)}
\mathrm{e}^{-(1+\nu')d(O,\gamma
O)} <\infty
\end{equation*}
for $\nu'>0$, 
 this completes the proof of the
convergence of the series
$\sum_{T} \dot\sigma_0(T)
C_0(\varphi, T)$.  This also
shows that $\sum_{T}\left|
\dot\sigma_0(T) C_0(\varphi,
T)\right|
\prec 
\left\Vert\varphi\right\Vert_\nu$,
which completes the proof of
Proposition~\ref{thm:SumCv}. 
\end{proof}

We will need at some point a
slightly stronger result:
\begin{lemma} 
\label{thm:SumCv:improv}
For $t$
sufficiently small, the series 
\begin{equation*}
\sum_{T} 
\dot\sigma_t(T)
\mathrm{e}^{\left|\alpha_t(T)\right|}
\mathrm{e}^{-(1+\nu)D_T}
\mathrm{e}^{-\left|u_T\right|}
\end{equation*}
converges, and this uniformly in
$t$.
\end{lemma}
\begin{proof} We again split the sum as 
\begin{equation*}
\sum_{T} 
\left\vert\dot\sigma_t(T)
\right\vert
\mathrm{e}^{\left|\alpha_t(T)
\right|}
\mathrm{e}^{-(1+\nu)D_T}
\mathrm{e}^{-\left|u_T\right|}
=
\sum_{\widetilde e}
\sum_{x_T\in \widetilde e}
\left\vert\dot\sigma_t(T)
\right\vert
\mathrm{e}^{\left|\alpha_t(T)
\right|}
\mathrm{e}^{-(1+\nu)D_T}
\mathrm{e}^{-\left|u_T\right|}.
\end{equation*}

If $x_T\in \widetilde e$, it
again follows from the
combinatorics of transverse
cocycles
\cite{Bon97b} that 
\begin{equation*}
\left|\alpha_t(T)\right|
\leq
\left|\alpha_t(T_{\tilde e})
\right| +
\left|\alpha_t(T)-
\alpha_t(T_{\tilde e})\right|
\prec 
\left\Vert\alpha_t\right\Vert_\Phi 
D_{\widetilde e} +
\left\Vert\alpha_t\right\Vert_\Phi 
r(T)
\end{equation*}
and
\begin{equation*}
\left|\dot\sigma_t(T)\right|
\leq
\left|\dot\sigma_t(T_{\tilde e})
\right| +
\left|\dot\sigma_t(T)-
\dot\sigma_t(T_{\tilde e})\right|
\prec 
\left\Vert
\dot\sigma_t\right\Vert_\Phi 
D_{\widetilde e} +
\left\Vert\dot
\sigma_t\right\Vert_\Phi  r(T).
\end{equation*}

Because $\alpha_0=0$ and by
existence of the derivative
cocycle
$\dot\alpha_0$, we have that
$\left\Vert\alpha_t
\right\Vert_\Phi \prec \left\vert
t\right\vert$. Also,
$\mathrm e^{-D_T} \asymp
\mathrm e^{-D_{\tilde
e}}$ and $\mathrm
e^{-\left\vert u_T\right
\vert } \prec \mathrm
e^{-at}$ for some $a>0$,
as in
Equations~\ref{equn:11.2}
and \ref{equn:11.4}. Therefore,
there are constants
$a$, $b>0$ such that 
\begin{equation*}
\begin{split}
\sum_{T} 
\left\vert\dot\sigma_t(T)
\right\vert
\mathrm{e}^{\left|\alpha_t(T)
\right|}
&\mathrm{e}^{-(1+\nu)D_T}
\mathrm{e}^{-\left|u_T\right|}\\
&\prec
\sum_{\widetilde e}
\sum_{x_T\in \widetilde e}
\left\| \dot\sigma_t\right\|_\Phi
\left( D_{\tilde e}+r(T)\right)
\mathrm{e}^{b\left\vert
t\right\vert 
D_{\tilde e}
+b\left\vert t\right\vert r(T)}
\mathrm{e}^{-(1+\nu)D_{\tilde e}}
\mathrm{e}^{-ar(T)}\\
&\prec 
\sum_{\widetilde e}
\sum_{x_T\in \widetilde e}
\mathrm{e}^{-(1+\nu')D_{\tilde e}}
\mathrm{e}^{-a'r(T)}
\end{split}
\end{equation*}
for $0<\nu'<\nu$ and $0<a'<a$, if
$t$ is chosen sufficiently small
that 
$b\left\vert
t\right\vert<\inf\{\nu-\nu',
a-a'\}$. We can then conclude the
proof as for
Proposition~\ref{thm:SumCv}.
\end{proof}

\section{Proof of the main
theorem}
The main theorem of the
article, namely
Theorem~\ref{thm:MainThm} in the
Introduction, is an immediate
consequence of
Theorems~\ref{thm:CompDeriv} and
\ref{thm:DerLin} below.
\begin{theorem}
\label{thm:CompDeriv}
Let $t\mapsto m_t$, $t\in\left]
-\varepsilon, \varepsilon\right[$,
be a differentiable curve in
$\mathcal T(S)$, let $\lambda$ be
a maximal geodesic lamination in
$S$, and let $\widetilde\lambda$
be its preimage in the universal
covering $\widetilde S$. Then,
for every H\"older
continuous function
$\varphi:G\bigl(\widetilde S\bigr)
\rightarrow \mathbb R$ with
compact support,
\begin{equation*}
\frac{d}{dt} \iint_{G(\widetilde
S)}
\varphi\thinspace dL_{m_t}\raise
-6pt\hbox{}_{\vert t=0} =
\sum_{T\not=T_O} \dot\sigma_0(T)
C_0(\varphi, T)
\end{equation*}
where: the sum is over all
components $T$ of $\widetilde
S-\widetilde\lambda$ which are
different from the component $T_O$
containing the base point
$O\in\widetilde S$; 
\begin{equation}
\label{eqn:DefC(T)bis}
C_0(\varphi, T) = \iint
_{G(\widetilde S)}
\varphi(h)
\left(
\cos\theta\left(g_3^T,h\right)
-\cos\theta\left(g_1^T,h\right)
-\cos\theta\left(g_2^T,h\right)
\right)
\thinspace dL_{m_0} (h)
\end{equation}
where the boundary of $T$
consists of the geodesics $g_1^T$,
$g_2^T$, $g_3^T$ with $g_3^T$
closest to $O$,
$\cos\theta\left(g,h\right)$ is
the cosine of the angle from the
geodesic $g$ to the geodesic $h$
at their intersection point,
measured counterclockwise for
the metric $m_0$, and
$\cos\theta\left(g,h\right)=0$ by
convention when $g$ and $h$ do not
meet;
$\dot\sigma_0(T)=
\frac{d}{dt}\sigma_t(T)
\raise -2pt\hbox{}_{|t=0}$ where
$\sigma_t(T)$ is the number
associated by the shearing
cocycle $\sigma_t\in\mathcal
H(\lambda)$ of $m_t\in\mathcal
T(S)$ to an arbitrary geodesic arc
joining the base point $O$ to $T$.
\end{theorem}

\begin{proof}[Proof of
Theorem~\ref{thm:CompDeriv}]
We are now ready to justify the
formal computations of
Section~\ref{sect:FormalComp}. We
will use the definitions and
notation of that section. Again,
we use the isometry $\Theta_{m_0}$
associated to $m_0$ to identify
$\widetilde S$ to
$\mathbb H^2$. In this
identification $\widetilde
S\cong\mathbb H^2$, the map
$\Theta_{m_0}$ is now the
identity and we use
the same symbol to represent an
object in
$\widetilde S$ and its image
under $\Theta_{m_0}$ in
$\mathbb H^2$.

Let $\mathcal U$ be a spanning
family of components of
$\widetilde S-\widetilde\lambda$,
and consider the product 
\begin{equation*}
E_{T_1T_2\dots T_n}^{\alpha_t}=
E_{T_1}^{\alpha_t(T_1)}
E_{T_2}^{\alpha_t(T_2)}
\dots
E_{T_n}^{\alpha_t(T_n)}
~
\prod_{U\in\mathcal U}
E_{g_3^U}^{\alpha_t(U)}
\end{equation*}
where $T_1$, $T_2$, \dots,
$T_n<\mathcal U$ as in
Equation~\ref{eqn:FiniteEarth},
and where $T_i<T_j\Rightarrow
i<j$.

Note that, for two
transverse cocycles $\alpha$,
$h$, 
\begin{equation*}
E_{T_1T_2\dots T_n}^{\alpha+h}=
E_{T_1'}^{h(T_1)}
E_{T_2'}^{h(T_2)}
\dots
E_{T_n'}^{h(T_n)}
~
\prod_{U\in\mathcal U}
E_{g_3^{U'}}^{h(U)}
E_{T_1T_2\dots T_n}^{\alpha}
\end{equation*}
where $T'_i=E_{T_1T_2\dots
T_n}^{\alpha}T_i$  and
$U'= E_{T_1T_2\dots
T_n}^{\alpha}U$. 
Since this product involves only
finitely many earthquakes, it
then follows from
Lemma~\ref{thm:ElemShear} that,
for every $t$,
\begin{multline}
\label{eqn:DerivFiniteShear}
\frac{d}{dt} \iint_{G(\mathbb
H^2)}
\varphi
\circ
\left(E_{T_1T_2\dots
T_n}^{\alpha_t}
\right)^{-1}
 dL_{\mathbb
H^2}
\\
=
\frac{d}{dt} \iint_{G(\mathbb
H^2)}
\varphi
\circ
\left(E_{T_1}^{\alpha_t(T_1)}
E_{T_2}^{\alpha_t(T_2)}
\dots
E_{T_n}^{\alpha_t(T_n)}
\prod_{U\in\mathcal U}
E_{g_3^U}^{\alpha_t(U)}
\right)^{-1}
 dL_{\mathbb
H^2}
\\
=
\sum_{i=1}^n
\dot\sigma_t(T_i)
C_0\bigl(
\varphi
\circ
\left(E_{T_1T_2\dots
T_n}^{\alpha_t}
\right)^{-1}, 
E_{T_1T_2\dots
T_n}^{\alpha_t}T_i\bigr)\\
+
\sum_{U\in\mathcal U}
\dot\sigma_t(U)
C_0\bigl(
\varphi
\circ
\left(E_{T_1T_2\dots
T_n}^{\alpha_t}
\right)^{-1}, 
E_{T_1T_2\dots
T_n}^{\alpha_t}g_3^U\bigr).
\end{multline}

As the finite family $\{T_1, T_2,
\dots, T_n\}$ converges to
$\{T\not=T_O; ,T<\mathcal U\}$, 
$E_{T_1T_2\dots
T_n}^{\alpha_t}$  converges to
\begin{equation*}
E_{\mathcal U}^{\alpha_t} =~
\overrightarrow
{\prod_{T<\mathcal U}}
E_T^{\alpha_t(T)}
~
\prod_{U\in\mathcal U}
E_{g_3^U}^{\alpha_t(U)}
\end{equation*}

We now pass to the limit in
Equation~\ref{eqn:DerivFiniteShear}
and conclude:
\begin{lemma}
\label{thm:DerComShear}
\begin{multline}
\label{eqn:DerivCompShear}
\frac{d}{dt} \iint_{G(\mathbb
H^2)}
\varphi\circ
\left(E_{\mathcal U}^{\alpha_t}
\right)^{-1}
 dL_{\mathbb
H^2}
\\=
\sum_{T<\mathcal U}
\dot\sigma_t(T)
C_0\bigl(
\varphi\circ
\left(E_{\mathcal U}^{\alpha_t}
\right)^{-1}, 
E_{\mathcal U}^{\alpha_t}T\bigr)+
\sum_{U\in\mathcal U}
\dot\sigma_t(U)
C_0\bigl(
\varphi\circ
\left(E_{\mathcal U}^{\alpha_t}
\right)^{-1},
E_{\mathcal
U}^{\alpha_t}g_3^U\bigr)
\end{multline}
for $t$ sufficiently small.
\end{lemma}
\begin{proof}[Proof of
Lemma~\ref{thm:DerComShear}] 
As $\{T_1, T_2,
\dots, T_n\}$ converges to
$\{T'\not=T_O; T'<\mathcal U\}$, 
the homeomorphism $E_{T_1T_2\dots
T_n}^{\alpha_t}$ of
$\partial\mathbb H^2$ uniformly
converges to
$E_{\mathcal U}^{\alpha_t}$. In
particular, for a fixed $T$, 
$C_0\bigl(
\varphi\circ
\left(E_{T_1T_2\dots
T_n}^{\alpha_t}
\right)^{-1},
E_{T_1 T_2,
\dots T_n}^{\alpha_t}T\bigr)$
converges to 
$C_0\bigl(
\varphi\circ
\left(E_{\mathcal U}^{\alpha_t}
\right)^{-1},
E_{\mathcal
U}^{\alpha_t}T\bigr)$. Similarly,
for
$U\in\mathcal U$,  
$C_0\bigl(
\varphi\circ
\left(E_{T_1T_2\dots
T_n}^{\alpha_t}
\right)^{-1},
E_{T_1 T_2,
\dots T_n}^{\alpha_t}g_3^U\bigr)$
converges to 
$C_0\bigl(
\varphi\circ
\left(E_{\mathcal U}^{\alpha_t}
\right)^{-1},
E_{\mathcal
U}^{\alpha_t}g_3^U\bigr)$. To
justify
Equation~\ref{eqn:DerivCompShear},
we need a uniform estimate for
the convergence of the sums in
the right hand sides of
Equations~\ref{eqn:DerivFiniteShear}
and \ref{eqn:DerivCompShear}.

For a component $T$ of
$\widetilde S-\widetilde
\lambda$, the refinement
Lemma~\ref{thm:EstC(T)bis} of
Proposition~\ref{thm:EstC(T)}
shows that 
\begin{equation}
\label{Eqn:EstC0FiniteEarth}
\begin{split}
 C_0\bigl(
\varphi\circ
\left(E_{T_1T_2\dots
T_n}^{\alpha_t}
\right)^{-1}, 
E_{T_1T_2\dots
T_n}^{\alpha_t}T\bigr)
&=
C_0\bigl(\varphi',
T'\bigr)\\
&\prec
\left\Vert\varphi'
\right\Vert_{g_3^{T'},
\nu}
\mathrm{e}^{-(1+\nu)D_{T'}}
\mathrm{e}^{-\left\vert
u_{T'}\right\vert}
\end{split}
\end{equation}
if we set $T'=E_{T_1T_2\dots
T_n}^{\alpha_t}T$ and
$\varphi'=\varphi\circ
\left(E_{T_1T_2\dots
T_n}^{\alpha_t}
\right)^{-1}$ to
simplify the notation.

\begin{sublemma} 
\label{thm:TandT'}
If $t$ is
sufficiently small, then
$
\left\vert D_{T'}-D_T\right\vert
\prec 1$ and $
\left\vert u_{T'}-u_T\right\vert
\prec
\left\vert\alpha_t(T)\right\vert+1
$.
\end{sublemma}
\begin{proof}[Proof of
Sublemma~\ref{thm:TandT'}]
 We distinguish cases,
according to whether $T$ is one
of the $T_i$ or not. 

If $T=T_{i_0}$, we can arrange
the indexing so that those $T_i$
with $T_i<T$ are exactly those
with $i<i_0$ (while keeping the
property that $T_i<T_j\Rightarrow
i<j$). Then,
\begin{equation*}
T'=I^{\alpha_t(T_1)}_{g_3^{T_1}}
I^{-\alpha_t(T_1)}_{g_{k_1}^{T_1}}
I^{\alpha_t(T_2)}_{g_3^{T_2}}
I^{-\alpha_t(T_2)}_{g_{k_2
}^{T_2}}
\dots
I^{\alpha_t(T_{i_0-1})
}_{g_3^{T_{i_0-1}}}
I^{-\alpha_t(T_{i_0-1}
)}_{g_{k_{i_0-1}}^{T_{i_0-1}}}
I^{\alpha_t(T_{i_0
})}_{g_3^{T_{i_0}}} T
\end{equation*}
where $g_{k_i}^{T_i}$ is the one
among the two geodesics
$g_{1}^{T_i}$ and
$g_{2}^{T_i}$ that separates $T$
from the base point $O$, and
where $I_g^a$ is the isometry of
$\mathbb H^2$ that acts by
translation of $a\in\mathbb R$
along the oriented geodesic $g$.

Let $T''=I^{\alpha_t(T_{i_0
})}_{g_3^{T_{i_0}}} T$. Then
$D_{T''}=D_T$ and $u_{T''}=u_T+
\alpha_t(T)$. On the other hand,
it is proved in \cite[\S
5]{Bon96} that, if $\alpha_t$ is
sufficiently small (depending on
the metric $m_0$) and therefore
if $t$ is sufficiently small,
\begin{equation*}
I^{\alpha_t(T_1)}_{g_3^{T_1}}
I^{-\alpha_t(T_1)}_{g_{k_1}^{T_1}}
I^{\alpha_t(T_2)}_{g_3^{T_2}}
I^{-\alpha_t(T_2)}_{g_{k_2}^{T_2}}
\dots
I^{\alpha_t(T_{i_0-1
})}_{g_3^{T_{i_0-1}}}
I^{-\alpha_t(T_{i_0-1
})}_{g_{k_{i_0-1}}^{T_{i_0-1}}}
\end{equation*}
remains bounded in the isometry
group of $\mathbb H^2$ (the bound
depends on
$\mathcal U$, but not the
condition on $t$).

Consequently, $\left\vert
D_{T'}-D_{T''}\right\vert$
and
$\left\vert
u_{T'}-u_{T''}\right\vert$
stay bounded, which proves
Sublemma~\ref{thm:TandT'} when $T$
is equal to some $T_{i_0}$.

The argument is similar when $T$
is equal to no $T_{i_0}$ (with
no term
$I^{\alpha_t(T_{i_0
})}_{g_3^{T_{i_0}}}
$), and leads to
$\left\vert
D_{T'}-D_{T}\right\vert \prec 1$
and $\left\vert
u_{T'}-u_{T}\right\vert \prec 1$
in this case.
\end{proof}

\begin{sublemma}
\label{thm:EstPhi'}
\begin{equation*}
\left\Vert \varphi'
\right\Vert_{g_3^{T'}, \nu}
\prec
\left\Vert \varphi
\right\Vert_\nu
\end{equation*}
\end{sublemma}
\begin{proof}[Proof of
Sublemma~\ref{thm:EstPhi'}]
By definition of
$\left\Vert\varphi'
\right\Vert_{g, \nu}$ in
Equation~\ref{Eqn:DefTNuNorm},
and since
$T'=E_{T_1T_2\dots
T_n}^{\alpha_t}T$ and
$\varphi'=\varphi\circ
\left(E_{T_1T_2\dots
T_n}^{\alpha_t}
\right)^{-1}$, 
\begin{equation*}
\left\Vert \varphi'
\right\Vert_{g_3^{T'}, \nu}
\leq
\left\Vert \varphi
\right\Vert_{g_3^T, \nu}
\max\bigl\{ 1, 
\mathrm e^{\nu\left( D_T -
D_{T'} \right)}\bigr \}
\prec \left\Vert \varphi
\right\Vert_\nu
\end{equation*}
where the second estimate comes
from Sublemma~\ref{thm:TandT'}
and from the fact that 
$\left\Vert \varphi
\right\Vert_{g_3^T, \nu}
\prec
\left\Vert \varphi
\right\Vert_\nu$. 
\end{proof}

Applying
Sublemmas~\ref{thm:TandT'} and
\ref{thm:EstPhi'} to Equation~\ref
{Eqn:EstC0FiniteEarth},
we conclude that 
\begin{equation}
\label{eqn:EstFinite Earth}
 C_0\bigl(
\varphi\circ
\left(E_{T_1T_2\dots
T_n}^{\alpha_t}
\right)^{-1}, 
E_{T_1T_2\dots
T_n}^{\alpha_t}T\bigr)
\prec
\left\Vert\varphi\right\Vert_\nu
\mathrm{e}^{\left\vert
\alpha_t(T)\right\vert}
\mathrm{e}^{-(1+\nu)D_{T}}
\mathrm{e}^{-\left\vert
u_{T}\right\vert}
\end{equation}

From
Equation~\ref{eqn:DerivFiniteShear},
we obtain 
\begin{multline}
\label{eqn:IntFiniteShear}
\iint_{G(\mathbb
H^2)}
\varphi
\circ
\left(E_{T_1T_2\dots
T_n}^{\alpha_t}
\right)^{-1}
 dL_{\mathbb
H^2}
-
\iint_{G(\mathbb
H^2)}
\varphi\,
 dL_{\mathbb
H^2}
\\=
\int_0^t
\sum_{i=1}^n
\dot\sigma_u(T_i)
C_0\bigl(\varphi\circ
\left(E_{T_1T_2\dots
T_n}^{\alpha_u}
\right)^{-1}, E_{T_1T_2\dots
T_n}^{\alpha_u}T_i\bigr)\,du\\
+
\int_0^t
\sum_{U\in\mathcal U}
\dot\sigma_u(U)
C_0\bigl(\varphi\circ
\left(E_{T_1T_2\dots
T_n}^{\alpha_u}
\right)^{-1}, E_{T_1T_2\dots
T_n}^{\alpha_u}g_3^U\bigr)
\,du.
\end{multline}

As $\{T_1, T_2,
\dots, T_n\}$ tends to
$\{T\not=T_O; ,T<\mathcal U\}$, 
the homeomorphism $E_{T_1T_2\dots
T_n}^{\alpha_t}$ of
$\partial\mathbb H^2$ uniformly
converges to
$E_{\mathcal U}^{\alpha_t}$. In
particular, for a fixed $T$, 
$C_0\bigl(
\varphi\circ
\left(E_{T_1T_2\dots
T_n}^{\alpha_u}
\right)^{-1},
E_{T_1 T_2,
\dots T_n}^{\alpha_u}T\bigr)$
converges to 
$C_0\bigl(
\varphi\circ
\left(E_{\mathcal U}^{\alpha_u}
\right)^{-1},
E_{\mathcal
U}^{\alpha_u}T\bigr)$ and 
$C_0\bigl(
\varphi\circ
\left(E_{T_1T_2\dots
T_n}^{\alpha_u}
\right)^{-1},
E_{T_1 T_2,
\dots T_n}^{\alpha_u}g_3^T\bigr)$
converges to 
$C_0\bigl(
\varphi\circ
\left(E_{\mathcal U}^{\alpha_u}
\right)^{-1},
E_{\mathcal
U}^{\alpha_u}g_3^T\bigr)$.
Combining
Equation~\ref{eqn:EstFinite
Earth} and
Lemma~\ref{thm:SumCv:improv}, we
conclude by dominated convergence
that
\begin{multline}
\label{eqn:IntCompShear}
 \iint_{G(\mathbb
H^2)}
\varphi\circ
\left(E_{\mathcal U}^{\alpha_t}
\right)^{-1}
 dL_{\mathbb
H^2} 
-
 \iint_{G(\mathbb
H^2)}
\varphi\,
 dL_{\mathbb
H^2} 
\\=
\int_0^t
\sum_{T<\mathcal U}
\dot\sigma_u(T)
C_0\bigl(\varphi\circ
\left(E_{\mathcal U}^{\alpha_u}
\right)^{-1},
E_{\mathcal U}^{\alpha_u}T\bigr)
\,du\\
+
\int_0^t
\sum_{U\in\mathcal U}
\dot\sigma_u(U)
C_0\bigl(\varphi\circ
\left(E_{\mathcal U}^{\alpha_u}
\right)^{-1},
E_{\mathcal
U}^{\alpha_u}g_3^U\bigr) \,du
\end{multline}
for $t$ sufficiently small.
Lemma~\ref{thm:DerComShear} then
follows by differentiating the
two sides of
Equation~\ref{eqn:IntCompShear}.
\end{proof}

We now let the spanning family
$\mathcal U$ uniformly tend to
infinity. Recall that this means
that the set $\{ T\not=T_O;
T<\mathcal U\}$ converges to the
set of all components $T\not=T_O$
of $\mathbb
H^2-\widetilde\lambda$. We will do
this in a controlled way.

Recall that $O_T$ denotes the
center of the ideal triangle $T$.
For $n\in\mathbb N$, let
$\mathcal V_n$ denote the set of
those components $T\not=T_O$ of
$G\bigl( \widetilde S\bigr)$ such
that $d(O,O_T)\leq n$. Note that
the set of all $O_T$ consists of
finitely many orbits of the
action of $\pi_1(S)$ on
$\widetilde S$, one orbit for each
component of $S-\lambda$, and is
consequently locally finite. It
follows that $\mathcal V_n$ is
finite. Let $\mathcal U_n \subset
\mathcal V_n$ consist of all
 elements of
$\mathcal V_n$ which are
maximal for~$<$. By construction,
$\mathcal U_n$ is a spanning
family.

As $n$ tends to $\infty$,
$\mathcal V_n$ converges to the
set of all components of
$\widetilde S-\widetilde\lambda$
which are different from $T_O$,
and it follows that $\mathcal
U_n$ uniformly converges to
infinity.

We will apply
Equation~\ref{eqn:DerivCompShear}
of Lemma~\ref{thm:DerComShear} to
$\mathcal U=\mathcal U_n$, and
let $n$ tend to $\infty$. 

\begin{lemma}
\label{thm:NoBdryTerm}
\begin{equation*}
\lim_{n\rightarrow \infty}
\sum_{U\in\mathcal U_n}
\dot\sigma_t(U)
C_0\bigl(\varphi\circ
\left(E_{\mathcal U_n}^{\alpha_t}
\right)^{-1}, E_{\mathcal
U_n}^{\alpha_t}g_3^U\bigr)=0
\end{equation*}
and this uniformly in $t$.
\end{lemma}
\begin{proof}[Proof of
Lemma~\ref{thm:NoBdryTerm}]
Let us begin with the case where
$t=0$. Using a train track $\Phi$
carrying $\lambda$ as in the
proof of
Proposition~\ref{thm:SumCv}, we
have that
$\dot\sigma_0(U)
\prec
\left\Vert\dot\sigma_0\right\Vert
_\Phi d(O,O_U)
\prec
\left\Vert\dot\sigma_0\right\Vert
_\Phi  n$ for every $U\in\mathcal
U_n$. 

By Proposition~\ref{thm:EstC(g)}, 
$C_0\bigl(\varphi,
g_3^U\bigr)
\prec 
\left\Vert
\varphi
\right\Vert_\nu
\mathrm{e}^{-(1+\nu)D_U}$.
The angle under which $g_3^U$ is
seen from the base point $O$ is
$\asymp
\mathrm{e}^{-D_U}$. Because the
 $U\in\mathcal U_n$ are pairwise
non comparable, the angular
sectors under which the
corresponding $g_3^U$ are seen
from
$O$ are disjoint. We conclude that
$\sum_{U\in\mathcal U_n} 
\mathrm{e}^{-D_U}
\prec 1$. Therefore,
\begin{multline*}
\sum_{U\in\mathcal U_n}
\dot\sigma_0(U)
C_0\bigl(\varphi, g_3^U\bigr)
\prec 
\left\Vert
\varphi
\right\Vert_\nu
\left\Vert
\dot\sigma_0\right\Vert_\Phi n
\sum_{U\in\mathcal U_n} 
\mathrm{e}^{-(1+\nu)D_U}\\
\prec 
\left\Vert
\varphi
\right\Vert_\nu
\left\Vert
\dot\sigma_0\right\Vert_\Phi n\,
\mathrm{e}^{-\nu\delta_n}
\sum_{U\in\mathcal U_n} 
\mathrm{e}^{-D_U}
\prec 
\left\Vert
\varphi
\right\Vert_\nu
\left\Vert
\dot\sigma_0\right\Vert_\Phi n\,
\mathrm{e}^{-\nu\delta_n}
\end{multline*}
if $\delta_n$ is the infimum of
the distance $D_U$ over
$U\in\mathcal U_n$.

If $\Delta$ is the diameter of
$S$ for the metric $m_0$, every
point of $\widetilde
S\cong\mathbb H^2$ is at distance
at most
$\Delta$ from some $O_T$. It
follows that, for every $T$,
there is a $T'>T$ such that
$d(O,O_{T'})\leq D_T + 2\Delta$.
Since $\mathcal U_n$ consists of
the maximal elements of $\mathcal
V_n$, we conclude that
$\delta_n\geq n-2\Delta$. As a
consequence, 
\begin{equation*}
\sum_{U\in\mathcal U_n}
\dot\sigma_0(U)
C_0\bigl(\varphi, g_3^U\bigr)
\prec 
\left\Vert
\varphi
\right\Vert_\nu
\left\Vert
\dot\sigma_0\right\Vert_\Phi n\,
\mathrm{e}^{-\nu n},
\end{equation*}
which tends to $0$ as $n$ tends to
$\infty$. 

For a general $t$, note that
$E_{\mathcal
U_n}^{\alpha_t}g_3^U=
E^{\alpha_t}g_3^U$ if
$U\in\mathcal U_n$. The
homeomorphism $E^{\alpha_t}
=
\Theta_{m_t}\circ
\Theta_{m_0}^{-1}$ of
$\partial_\infty \mathbb H^2$ is
H\"older continuous with
uniformly bounded H\"older
exponent and norm (depending on a
compact subset of Teichm\"uller
space containing the $m_t$).
Since the angle under which 
$E^{\alpha_t}g_3^U$ is seen from
$O$ is $\asymp
\mathrm{e}^{-D_U^t}$, where
$D_U^t$ is the distance from $O$
to $E^{\alpha_t}g_3^U$, it
follows that $D_U^t\geq \mu D_U$
for some constant $\mu>0$.
In particular, the infimum
$\delta_n^t$ of all $D_U^t$ with
$U\in\mathcal U_n$ is such that 
$\delta_n^t\geq\mu\delta_n \geq
\mu n - 2\mu \Delta$.
Proposition~\ref{thm:EstC(g)}
gives that
\begin{equation*}
C_0\left( \varphi\circ
\bigl(E_{\mathcal U_n}^{\alpha_t}
\right)^{-1}, 
E_{\mathcal
U_n}^{\alpha_t} g_3^U\bigr)
\prec
\bigl\Vert
\varphi\circ
\left(E_{\mathcal U_n}^{\alpha_t}
\right)^{-1}
\bigr\Vert_{
E_{\mathcal
U_n}^{\alpha_t}g_3^U, \nu}
\mathrm e^{-(1+\nu)D_U^t}.
\end{equation*}
 By
definition of $\left\Vert
\ \right\Vert_{g,
\nu}$, we have 
\begin{equation*}
\bigl\Vert
\varphi\circ
\left(E_{\mathcal U_n}^{\alpha_t}
\right)^{-1}
\bigr\Vert_{
E_{\mathcal
U_n}^{\alpha_t}g_3^U, \nu}
\leq
\bigl\Vert
\varphi
\bigr\Vert_{
g_3^U, \nu}
\max \bigl\{ 1, \mathrm
e^{\nu(D_U^t-D_U)} \bigr\}
\prec
\bigl\Vert
\varphi
\bigr\Vert_{ \nu}
\max \bigl\{ 1, \mathrm
e^{\nu(D_U^t-D_U)} \bigr\}.
\end{equation*}
Decomposing $\mathcal U_n$ as the
union of the set $\mathcal U_n^+$
of those $U\in \mathcal U_n$ with
$D_U^t\geq D_U$ and the set
$\mathcal U_n^-$
of those $U$ with
$D_U^t < D_U$, we conclude that
\begin{multline*}
\sum_{U\in\mathcal U_n}
\dot\sigma_t(U)
C_0\bigl(\varphi\circ
\left(E_{\mathcal U_n}^{\alpha_t}
\right)^{-1}, E_{\mathcal
U_n}^{\alpha_t}g_3^U\bigr)\\
\prec
\left\Vert \varphi
\right\Vert_\nu 
\left\Vert\dot\sigma_0\right\Vert
n
\sum_{U\in\mathcal U_n^+} 
\mathrm{e}^{-D_U^t-\nu D_U}
+
\left\Vert \varphi
\right\Vert_\nu 
\left\Vert\dot\sigma_0\right\Vert
n
\sum_{U\in\mathcal U_n^-} 
\mathrm{e}^{-D_U^t-\nu D_U^t}
\\
\prec
\left\Vert \varphi
\right\Vert_\nu 
\left\Vert\dot\sigma_0\right\Vert
n\,
\bigl(
\mathrm{e}^{-\nu\delta_n}
+
\mathrm{e}^{-\nu\delta_n^t}
\bigr)
\sum_{U\in\mathcal U_n} 
\mathrm{e}^{-D_U^t}\\
\prec
\left\Vert \varphi
\right\Vert_\nu 
\left\Vert\dot\sigma_0\right\Vert
n\,
\bigl(
\mathrm{e}^{-\nu n}
+
\mathrm{e}^{-\nu\mu n}
\bigr)
\end{multline*}
since $\sum_{U\in\mathcal U_n} 
\mathrm{e}^{-D_U^t}
\prec 1$, as in the case where
$t=0$. This completes the proof
of
Lemma~\ref{thm:NoBdryTerm}.
\end{proof}

\begin{lemma}
\label{thm:CvPartialEarth}
For $t$ sufficiently small, the
sum  
\begin{equation*}
\sum_T
\dot\sigma_t(T)
C_0\bigl(\varphi\circ
\left(E_{\mathcal U}^{\alpha_t}
\right)^{-1}, E_{\mathcal
U}^{\alpha_t}T\bigr)
\end{equation*} converges,
uniformly in $\mathcal U$ and in
$t$.
\end{lemma}
\begin{proof} By
Lemma~\ref{thm:EstC(T)bis} and by
definition of $\left\Vert
\phantom{m} \right\Vert_{g, \nu}$,
\begin{equation*}
\begin{split}
C_0\bigl(\varphi\circ
\left(E_{\mathcal U}^{\alpha_t}
\right)^{-1}, E_{\mathcal
U}^{\alpha_t}T\bigr)
&\prec
\bigl\Vert 
\varphi\circ
\left(E_{\mathcal U}^{\alpha_t}
\right)^{-1}
\bigr\Vert
_{g_3^{T'}, \nu}
\mathrm e^{-(1+\nu)D_{T'}}
\mathrm e^{-\left\vert
u_{T'}\right\vert} \\
&\prec
\bigl\Vert 
\varphi
\bigr\Vert
_{\nu}
\max\bigr\{1, \mathrm
e^{\nu(D_{T'}-D_T)}
\bigr\}
\mathrm e^{-(1+\nu)D_{T'}}
\mathrm e^{-\left\vert
u_{T'}\right\vert} 
\end{split}
\end{equation*} if we set
$T'=E_{\mathcal U}^{\alpha_t}T$.

Note that $T' = E^{\alpha_t}T$
when 
$T<\mathcal U$, and $T'=T$
otherwise. In particular, in the
first case, $D_{T'}$ is the
$m_t$--distance from the base
point $O$ to the $m_t$--geodesic
triangle in $\widetilde S$
corresponding to $T$. As a
consequence,
$D_{T'}
\asymp D_T$ in both cases, and we
conclude that 
\begin{equation*}
C_0\bigl(\varphi\circ
\left(E_{\mathcal U}^{\alpha_t}
\right)^{-1}, E_{\mathcal
U}^{\alpha_t}T\bigr)
\prec
\bigl\Vert 
\varphi
\bigr\Vert
_{\nu}
\mathrm e^{-(1+\nu')D_{T'}}
\mathrm e^{-\left\vert
u_{T'}\right\vert} 
\end{equation*}
for some $\nu'>0$ depending on
$\nu$ and on  a compact subset
of $\mathcal T(S)$ containing the
image of the curve
$t\mapsto m_t$. Therefore, 
$\sum_T
\dot\sigma_t(T)
C_0\bigl(\varphi\circ
\left(E_{\mathcal U}^{\alpha_t}
\right)^{-1}, E_{\mathcal
U}^{\alpha_t}T\bigr)$ is bounded
by the sum of
$\bigl\Vert 
\varphi
\bigr\Vert
_{\nu} \sum_T
\dot\sigma_t(T)
\mathrm e^{-(1+\nu')D_{T}}
\mathrm e^{-\left\vert
u_{T}\right\vert}$ and 
$\bigl\Vert 
\varphi
\bigr\Vert
_{\nu} \sum_{T}
\dot\sigma_t(T)
\mathrm e^{-(1+\nu')D_{T''}}
\mathrm e^{-\left\vert
u_{T''}\right\vert}$, where $T$
ranges over all components of
$\widetilde S -\widetilde\lambda$
that are different from $T_O$, 
and where
$T'' = E^{\alpha_t}T$ in the
second sum. These two sums
converge by
Lemma~\ref{thm:SumCv:improv},
applied to the metric $m_t$
instead of
$m_0$ for the second sum. This
proves
Lemma~\ref{thm:CvPartialEarth}.
\end{proof}

If we apply
Equation~\ref{eqn:IntCompShear}
to $\mathcal U=\mathcal U_n$ and
let $n$ tend to $\infty$ we 
then conclude, by dominated
convergence  and by
Lemma~\ref{thm:NoBdryTerm}, that
\begin{multline}
\label{eqn:IntFullShear}
 \iint_{G(\mathbb
H^2)}
\varphi\circ
\left(E^{\alpha_t}
\right)^{-1}
 dL_{\mathbb
H^2} 
-
 \iint_{G(\mathbb
H^2)}
\varphi\,
 dL_{\mathbb
H^2} 
\\=
\int_0^t
\sum_{T}
\dot\sigma_u(T)
C_0\bigl(\varphi\circ
\left(E^{\alpha_u}
\right)^{-1}, E^{\alpha_u}T\bigr)
\,du .
\end{multline}

Differentiating both sides of
Equation~\ref{eqn:IntFullShear},
it follows that
\begin{multline*}
\frac{d}{dt} \iint_{G(\widetilde
S)}
\varphi\thinspace dL_{m_t}\raise
-6pt\hbox{}_{\vert t=0} \\
=
\frac{d}{dt} \iint_{G(\mathbb
H^2)}
\varphi\circ
\left(E^{\alpha_t}
\right)^{-1}
 dL_{\mathbb
H^2} \raise
-6pt\hbox{}_{\vert t=0} =
\sum_{T} \dot\sigma_0(T)
C_0(\varphi, T)
\end{multline*}
which completes the proof of
Theorem~\ref{thm:CompDeriv}.
\end{proof}

\begin{theorem}
\label{thm:DerLin}
With the data of
Theorem~\ref{thm:CompDeriv}, the
map $T_{m_0}L(\dot m_0)$ which to
a H\"older continuous function
$\varphi:G\bigl(\widetilde S\bigr)
\rightarrow \mathbb R$ with
compact support associates 
\begin{equation*}
T_{m_0}L(\dot m_0)(\varphi)=
\frac{d}{dt} \iint_{G(\widetilde
S)}
\varphi\thinspace dL_{m_t}\raise
-6pt\hbox{}_{\vert t=0} 
\end{equation*}
is a H\"older geodesic current
$T_{m_0}L(\dot m_0) \in\mathcal
H(S)$ which depends only on the
tangent vector $\dot m_0 =
\frac{d}{dt} m_t\raise
-2pt\hbox{}_{|t=0} \in
T_{m_0}\mathcal T(S)$. The map
\begin{equation*}
T_{m_0}L :
T_{m_0}\mathcal T(S) \rightarrow
\mathcal H(S)
\end{equation*}
so defined is linear, and depends
continuously on $m_0$. 
\end{theorem}

\begin{proof}
By Theorem~\ref{thm:CompDeriv}
and Proposition~\ref{thm:SumCv},
$T_{m_0}L(\dot m_0)$ defines a
H\"older distribution on the
space $G\bigl(\widetilde
S\bigr)$. Also, since $L_{m_t}$
is invariant under the action of
$\pi_1(S)$, so is its derivative
$T_{m_0}L(\dot m_0)$. (Note that
it is less obvious to see the
invariance directly from the
formula $T_{m_0}L(\dot
m_0)(\varphi)=\sum_{T}
\dot\sigma_0(T) C_0(\varphi, T)$).
As a consequence, $T_{m_0}L(\dot
m_0)$ is a H\"older geodesic
current.

Since the map which to a
hyperbolic metric $m\in\mathcal
T(S)$ associates its shearing
cocycle $\sigma\in\mathcal
H(\lambda)$ is a local
diffeomorphism,
$\dot\sigma_0\in\mathcal
H(\lambda)$ depends linearly
on the tangent vector $\dot m_0
\in T_{m_0}\mathcal T(S)$. It
follows that $T_{m_0}L(\dot
m_0)(\varphi)=\sum_{T}
\dot\sigma_0(T) C_0(\varphi, T)$
depends linearly on $\dot m_0$.
In other words, the map
$T_{m_0}L$ is linear.

Finally, for two geodesics $g$,
$h\in G\bigl(\widetilde S\bigr)$,
the cosine $\cos\theta(g,h)$
depends continuously on the
metric $m_0$ used to measure the
angle $\theta(g,h)$. By bounded
convergence, it follows that each
$C_0(\varphi, T)$ is a continuous
function of $m_0$. In addition,
the arguments of the proof of
Proposition~\ref{thm:SumCv} show
that the convergence of the series
$\sum_{T}
\dot\sigma_0(T) C_0(\varphi, T)$
is uniform as $m_0$ stays in a
compact subset of $\mathcal
T(S)$. Consequently, for every
$\varphi$,  $T_{m_0}L(\dot
m_0)(\varphi)=\sum_{T}
\dot\sigma_0(T) C_0(\varphi, T)$
depends continuously on $m_0$. In
other words, $T_{m_0}L(\dot
m_0)\in\mathcal H(S)$ is a
continuous function of $m_0$ and
$\dot m_0$. 
\end{proof}


\begin{thebibliography}{CEG}

\bibitem[BS]{BirSer} Joan S.
Birman, Caroline Series,
\emph{Geodesics with bounded
intersection number on surfaces
are sparsely distributed},
Topology \textbf{24} (1985),
217--225.

\bibitem[Bo1]{Bon88}
Francis Bonahon,  \emph{The
geometry of {T}eichm\"uller space
via geodesic
currents}, Invent. Math.
\textbf{92} (1988), 139--162. 

\bibitem[Bo2]{Bon96}
\bysame, \emph{Shearing hyperbolic 
surfaces, bending pleated
surfaces and {T}hurston's
symplectic form}, Ann. Fac. Sci.
Toulouse Math. (6)
\textbf{5}
(1996),  233--297.

\bibitem[Bo3]{Bon97a}
\bysame,
\emph{Geodesic laminations with
transverse {H}\"older
distributions}, Ann. Sci. \'Ecole
Norm. Sup. (4)
\textbf{30} (1997),  205--240.
  

\bibitem[Bo4]{Bon97b}
\bysame,
\emph{Transverse {H}\"older
distributions for geodesic
laminations}, Topology
\textbf{36} (1997),  103--122.

\bibitem [CEG]{CanEpsGre}
Richard D. Canary, David B. A.
Epstein, Paul Green, \emph{Notes
on notes of Thurston}, in:
\emph{Analytical and geometric
aspects of hyperbolic space
(Coventry/Durham, 1984)}, 3--92,
London Math. Soc. Lecture Note
Ser. vol. \textbf{111}, Cambridge
Univ. Press, Cambridge, 1987.

\bibitem[Gr]{Gro}
Mikhael L. Gromov,
\emph{Hyperbolic groups}, in:
\emph{Essays in group theory},
75--263, Math. Sci. Res. Inst.
Publ. vol. \textbf{8}, Springer,
New York, 1987.

\bibitem[KH]{KatHas95} Anatole
Katok, Boris Hasselblatt,
\emph{Introduction to the modern
theory of dynamical systems},
Cambridge University Press,
Cambridge, 1995.

\bibitem[Ot]{Ota90}
Jean-Pierre Otal, \emph{Le
spectre marqu\'e des longueurs des
surfaces \`a courbure n\'egative},
Ann. of Math. \textbf{131} (1990),
151--162. 

\bibitem[PH]{PenHar}
Robert C. Penner, John L. Harer,
\emph{Combinatorics of train
tracks}, Annals of Mathematics
Studies vol. \textbf{125},
Princeton University Press,
Princeton, 1992. 

\bibitem[\v Sa] {Sar} Dragomir 
\v Sari\'c, \emph{Infinitesimal
Liouville distributions for
Teichm\"uller space}, preprint,
University of Southern
California, 2002. 

\bibitem[Su]{Sul76}
Dennis P. Sullivan, \emph{Cycles
for the dynamical study of
foliated manifolds and complex
manifolds}, Invent. Math.
\textbf{36} (1976), 225--255.

\bibitem[Th]{Thu86d}
William~P. Thurston, 
\emph{Minimal stretch maps
between hyperbolic surfaces},
  Unpublished preprint.

\bibitem[Wo]{Wol} Scott A.
Wolpert,
\emph{ Thurston's 
Riemannian metric for 
Teichm{\"u}ller space},  J. 
Differ. Geom. \textbf{23} (1986),
143--174.

\end{thebibliography}

\end{document}